\documentclass[a4paper, 11pt]{article}

\pdfoutput=1

\usepackage{color}
\usepackage{graphicx}
\usepackage{amsmath}
\usepackage{amsfonts}
\usepackage{amssymb}
\usepackage{amsthm}
\usepackage{epstopdf}
\usepackage{algorithmic}
\usepackage{setspace}
\usepackage[cm]{fullpage}
\usepackage{multirow}

\theoremstyle{remark}
\newtheorem{remark}{Remark}
\newtheorem{algorithm}{Algorithm}
\newtheorem{problem}{Problem}

\newcommand{\vlambda}{{\boldsymbol{\lambda}}}

\newcommand{\dabR}{\mathbb{R}}

\newcommand{\Image}{\textit{Im\,}}
\newcommand{\Ker}{\textit{Ker\,}}
\newcommand{\diag}{\textit{diag\,}}

\newcommand{\mbf}[1]{\mbox{\boldmath$#1$}}

\title{A TFETI Domain Decomposition Solver for Elastoplastic Problems}

\author{M. \v{C}ERM{\'A}K$^1$, T. KOZUBEK$^1$, S. SYSALA$^2$, J. VALDMAN$^{1,3}$\footnote{corresponding author, email: jan.valdman@vsb.cz }
\\ \\
Centre of Excellence IT4Innovations, \\
$^1$V\v SB-Technical University of Ostrava, \\
$^2$Institute of Geonics AS CR, v.v.i., Ostrava, \\
$^3$The Institute of Information Theory and Automation, Prague
}

\begin{document}
\maketitle

\begin{abstract}
We propose an algorithm for the efficient parallel implementation of elastoplastic problems with hardening based on the so-called TFETI (Total Finite Element Tearing and Interconnecting) domain decomposition method. We consider an associated elastoplastic model with the von Mises plastic criterion and the linear isotropic hardening law. Such a model is discretized by the implicit Euler method in time and the consequent one time step elastoplastic problem by the finite element method in space. The latter results in a system of nonlinear equations with a strongly semismooth and strongly monotone operator. The semismooth Newton method is applied to solve this nonlinear system. Corresponding linearized problems arising in the Newton iterations are solved in parallel by the above mentioned TFETI domain decomposition method. The proposed TFETI based algorithm was implemented in Matlab parallel environment and its performance was illustrated on a 3D elastoplastic benchmark. Numerical results for different time discretizations and mesh levels are presented and discussed and a local quadratic convergence of the semismooth Newton method is observed.
\end{abstract}

\section{Introduction}
Elastoplastic processes describe behaviour of solid continuum beyond reversible elastic deformations. They are typically described by hysteresis models with a time memory \cite{BS96, K96}.
The rigorous mathematical analysis of elastoplastic problems and the numerical methods for their solution started to appear in the late 70ies and in the early 80ies by the work of C. Johnson \cite{J76, J78}, H. Matthies \cite{Ma79, Ma79b}, V. Korneev and U. Langer \cite{KL84}, J. Ne\v{c}as and I. Hlav\'{a}\v{c}ek \cite{NeHl81} and others. Since then a lot of mathematical contributions to computational plasticity have been written, we refer at least to the monographs by J. Simo and T. Hughes \cite{SimoHughes} and W. Han and B. Reddy \cite{HaRe99}.

In this paper, we focus on the efficient parallel implementation of elastoplastic problems based on the TFETI domain decomposition method \cite{TFETI2006,KVMHDHKC}. More specifically, we consider an associated elastoplasticity with the von Mises plastic criterion and the linear isotropic hardening law (see e.g. \cite{HaRe99,BLA99,NPO08}). The corresponding elastoplastic constitutive model is discretized by the implicit Euler method in time and consequently a nonlinear stress-strain relation is implemented by the return mapping concept (see e.g. \cite{NPO08, ACZ99, BLA99}). This approach together with the balance equation, the small strain assumption and a combination of the Dirichlet and Neumann boundary conditions leads to the solution of a nonlinear variational equation with respect to the primal unknown displacement in each time step. Such an equation can also be equivalently formulated as a minimization problem with a potential energy functional (see e.g. \cite{GV09, Sy09}).

By a finite element space discretization of the one time step problem, we obtain a system of nonlinear equations. The corresponding nonlinear operator is nondifferentiable but strongly semismooth. Therefore, it is suitable to choose the semismooth Newton method for solving the system since the strong semismoothness together with other properties ensure local quadratic convergence. Semismooth functions in finite dimensional spaces and the semismooth Newton method were introduced in \cite{QS93}. In elastoplasticity, the semismoothness was investigated for example in \cite{GV09, SW11, Sy09, Sy12}.

In each Newton iteration, it is necessary to solve the respective linearized problem. Different linear solvers including those based on multigrid have been successfully tested in \cite{Wi10, GKLSV12}. Moreover, since the linear systems of equations corresponding to the elastic and elastoplastic problems are spectrally equivalent \cite{KLV04}, all preconditioners for elastic problems can be applied to elastoplastic ones as well.

A linear solver considered in this paper is based on a FETI type domain decomposition method enabling its efficient parallel implementation. The standard FETI method (FETI-1) was originally introduced by Farhat and Roux \cite{fr92} and theoretically analyzed by Mandel and Tezaur \cite{ManTez-1996}. A systematic overview of 
 Using this approach, a body is partitioned into non-overlapping subdomains, an elliptic problem with Neumann boundary conditions is defined
for each subdomain, and intersubdomain field continuity is enforced via Lagrange multipliers.
The Lagrange multipliers are efficiently solved from a dual problem by a variant of the conjugate gradient algorithm.
The first practical implementations exploited only the favorable distribution of the spectrum of the matrix of the smaller problem \cite{r92}, known also
as the dual Schur complement matrix, but such algorithm was efficient only with a small number of subdomains. Later, Farhat, Mandel, and Roux introduced
a ``natural coarse problem'' whose solution was implemented by auxiliary projectors so that the resulting algorithm became in a sense optimal \cite{fmr94,r98}. Here, we use the Total-FETI (TFETI) \cite{TFETI2006,KVMHDHKC} variant of FETI domain decomposition method, where even the Dirichlet boundary conditions are enforced by Lagrange multipliers. Hence all subdomain stiffness matrices are singular with a-priori known kernels which is a great advantage in the numerical solution. With known kernel basis we can regularize effectively the stiffness matrix without extra fill in and use any standard sparse Cholesky type decomposition method for nonsingular matrices \cite{KucKozMar-2013,BrzDosKovKozMar-2010}.

 For systematic overview of domain decomposition methods we refer to \cite{TosWid-2005}. We mention  FETI-DP and BDDC methods as examples of  alternatives to the Total-FETI method. The FETI-DP (dual-primal FETI) is a popular
way of avoiding singular stiffness matrices of the local problems and was first introduced by Farhat, Lesoinne, Le Tallec, Pierson, and Rixen \cite{FLLPR01}. On the other hand, in \cite{KucKozMar-2013} we gave numerical examples showing that the matrices of the systems arising in TFETI may be better conditioned than those arising in FETI-DP. More theoretical and implementation details about FETI-DP and its extensions and improvements are given in \cite{FarLesPie-2000,KlaWid-2001,KlaRhe-2006}. The BDDC (balancing domain decomposition by constraints) was introduced by Dohrmann \cite{Dohrmann-2003} as a simpler primal alternative to the FETI-DP. The name of the method was coined by Mandel and Dohrmann \cite{ManDoh-2003}, because it can be understood as further development of the BDD (balancing domain decomposition) method \cite{Mandel-1993}. An alternative approach to the solution of elastoplastic problems based on the Schwarz domain decomposition method with overlap was introduced in \cite{BG94}.

The structure of the paper is as follows: In Section 2, we introduce a quasistatic scheme of a solid mechanics problem. Within this context, we consider both elastic and elastoplastic models. The elastic model is introduced for methodical purposes because it is an essential part of the investigated elastoplastic model. After its time discretization we summarize the resulting one time step problem. In Section 3, the finite element space discretization of the one time step problem and its nonlinear algebraic formulation are described in details. The semismooth Newton method is applied to treat this nonlinearity. The TFETI method combined with the projected conjugate gradient algorithm is derived in Section 4 to solve the linearized problems appearing in the Newton iterations. Finally, the algorithm for solving the whole elastoplastic problem is summarized. In Section 5, the performance of the proposed algorithm is illustrated on numerical experiments. Final comments are summarized in Section 6.

\section{Elastic and elastoplastic models} \label{section.elastoplastic}

In this section, we summarize elastic and elastoplastic models in a quasistatic framework. Firstly, we introduce basic notation and assumptions that will be used for setting the models. Secondly, we describe the linear elastic constitutive model for an isotropic material. Thirdly, we introduce the elastoplastic initial value constitutive model based on the described elasticity, the von Mises yield function, the associated plastic flow rule and the linear isotropic hardening represented by the accumulated plastic strain. Fourthly, we consider time discretization of the constitutive model given by the implicit Euler method. The explicit form of the time discretized elastoplastic operator is introduced.

\subsection{Notation and assumptions}

Let us consider a deformable body occupying a domain $\Omega\subset \mathbb R^3$ with a Lipschitz continuous boundary $\Gamma=\partial\Omega$. We will describe the state of the body during a loading process by the Cauchy stress tensor $\sigma\in S$, the displacement $u\in\mathbb R^3$ and the small strain tensor $\varepsilon\in S$. Here $S=\mathbb R^{3\times 3}_{sym}$ is the space of all symmetric second order tensors. Other variables that are necessary for defining the elastoplastic models will be introduced in Subsection \ref{subsection_plast_model}. More details can be found in \cite{NPO08}.

The above variables depend on the spatial variable $x\in\Omega$ and on the time variable $t\in Q=[t_0,t^*]$. The small strain tensor is related to the displacement by the linear relation
\begin{equation}
\varepsilon(u)=\frac{1}{2}\left(\bigtriangledown u +(\bigtriangledown u)^T\right).
\label{strain-displ}
\end{equation}
The equilibrium equation in the quasistatic case reads
\begin{equation}
-\mbox{div}(\sigma(x,t))=g(x,t)\quad\forall (x,t)\in\Omega\times Q,
\label{balance}
\end{equation}
where $g(x,t)\in\mathbb R^3$ represents the volume force acting at the point $x\in\Omega$ and the time $t\in Q$.

\begin{figure}[n]   
\begin{center}
\includegraphics[height=5cm]{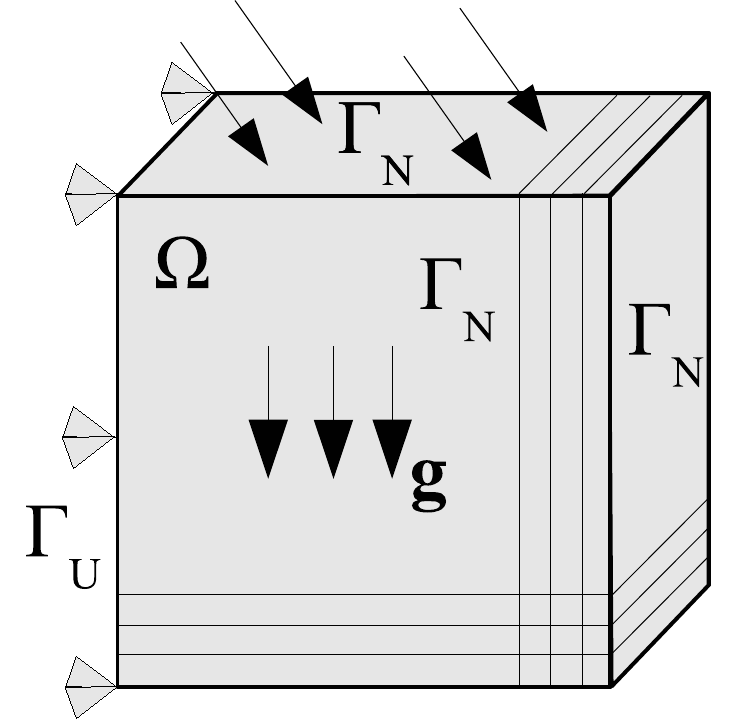}
  \caption{\label{fig1}Geometry of the domain $\Omega$ with imposed boundary conditions.}
\end{center}
\end{figure}

Let the boundary $\Gamma$ be fixed on a part $\Gamma_U$ that has a nonzero Lebesgue measure with respect to $\Gamma$, i.e., we prescribe the homogeneous Dirichlet boundary condition on $\Gamma_U$:
\begin{equation}
u(x,t)=0\quad\forall (x,t)\in\Gamma_U\times Q.
\label{dirichlet}
\end{equation}
On the rest of the boundary $\Gamma_N=\Gamma\setminus\Gamma_U$, we prescribe the Neumann boundary conditions
\begin{equation}
\sigma(x,t)n(x)=F(x,t)\quad\forall (x,t)\in\Gamma_N\times Q,
\label{neumann}
\end{equation}
where $n(x)$ denotes the exterior unit normal and $F(x,t)$ denotes a prescribed surface forces at the point $x\in\Gamma_N$ and the time $t\in Q$. Geometry of $\Omega$ with imposed boundary conditions is depicted in Figure \ref{fig1}. Similarly, we can consider other boundary conditions, for example symmetry and periodic conditions.

For a weak formulation of the investigated problems, it is sufficient to introduce the space of kinematically admissible displacements,
\begin{equation}
V=\left\{v\in[H^1(\Omega)]^3:\;v=0\;\; \mbox{on}\;\;\Gamma_U\right\}.
\label{V-space}
\end{equation}
Then the conditions (\ref{balance})--(\ref{neumann}) can be written in a weak sense by
\begin{equation}
\int_\Omega\langle\sigma,\varepsilon(v)\rangle_F dx=\int_\Omega g^Tv dx + \int_{\Gamma_N}F^Tv ds\quad\forall v\in V,\;\forall t\in Q.
\label{weak-balance}
\end{equation}
Here $\varepsilon(v)$ is defined by (\ref{strain-displ}), $\langle.,.\rangle_F$ and $\|.\|_F$ denote the Frobenius scalar product and the corresponding norm on the space $S$, respectively. We assume that the functions $\sigma,F,g$ are sufficiently smooth such that the integrals in (\ref{weak-balance}) are correctly defined in the Lebesgue sense.

To complete the investigated (generally quasistatic) models, we will prescribe the constitutive relations between the stress, the strain and eventually other variables, see Subsection \ref{subsection_elast_model} and  \ref{subsection_plast_model}.

\subsection{Elastic model}
\label{subsection_elast_model}

We consider the elastic constitutive model given by the Hooke law for isotropic material,
\begin{equation}
\sigma=\mathbb C\varepsilon=\lambda\mbox{tr}(\varepsilon)I+2\mu\varepsilon
\label{Hook_law}
\end{equation}
with the Lame coefficients $\lambda, \mu$. For the sake of simplicity, we assume a homogeneous material, i.e., the constant coefficients $\lambda, \mu>0$. The trace operator of a tensor is denoted by $\mbox{tr}(.)$ and $I$ denotes the identity.

It will be useful to introduce the volumetric and deviatoric parts of a tensor $\eta\in S$ by
\begin{equation}
\mbox{vol}(\eta)=\frac{1}{3}\mbox{tr}(\eta)I,\quad \mbox{dev}(\eta)=\eta-\mbox{vol}(\eta).
\end{equation}
It holds that
\begin{equation}
\langle\mbox{vol}(\eta),\mbox{dev}(\xi)\rangle_F=0,\;\;\langle\mbox{dev}(\eta),\xi\rangle_F=\langle\mbox{dev}(\eta),\mbox{dev}(\xi)\rangle_F\quad\forall\eta,\xi\in S.
\label{vol-dev_prop}
\end{equation}
By (\ref{vol-dev_prop}), we can find that the fourth order tensor $\mathbb C$, defined by (\ref{Hook_law}), is symmetric and elliptic, i.e.,
\begin{equation}
\langle \mathbb C\eta,\xi\rangle_F=\langle \eta,\mathbb C\xi\rangle_F,\;\;\langle \mathbb C\eta,\eta\rangle_F\geq 2\mu\|\eta\|_F^2\quad\forall\eta,\xi\in S.
\label{C-prop}
\end{equation}

If we substitute (\ref{Hook_law}) into (\ref{weak-balance}), we obtain for any fixed $t\in Q$ the weak formulation of the elastic problem.
\begin{problem}[Elastic problem]{\it Find $u=u(x,t) \in V$ such that
\begin{equation}
a_e(u,v)=\int_\Omega g^Tv dx + \int_{\Gamma_N}F^Tv ds\quad\forall v\in V,
\label{weak-elast}
\end{equation}
where the bilinear form on $V$ reads
\begin{equation}
a_e(w,v)=\int_\Omega\langle\mathbb C\varepsilon(w),\varepsilon(v)\rangle_F dx,\quad  w,v\in V.
\end{equation}
}
\end{problem}
Due to (\ref{C-prop}) and the Korn inequality \cite{NeHl81}, $a_e(w,v)$ is symmetric and $V$-elliptic
\begin{equation}
a_e(w,v)=a_e(v,w),\quad\exists c>0:\;\;a_e(v,v)\geq c\|v\|_V^2\quad\forall v,w\in V.
\label{a_e-est}
\end{equation}
The mentioned properties of the form $a_e$ ensure that the elastic problem (\ref{weak-elast}) has a unique solution $u\in V$, for example by Lax-Milgram lemma \cite{NeHl81}. Notice that the problem (\ref{weak-elast}) does not depend on the load history, so it is a static problem.

\subsection{Elastoplastic initial value constitutive model}
\label{subsection_plast_model}

In comparison to elasticity (see \eqref{Hook_law}), elastoplasticity is a time-dependent model where the history of loading is taken into account. We will assume
associated elastoplasticity with von Mises plastic criterion and linear isotropic hardening law (see e.g. \cite{HaRe99, BLA99, NPO08}). For details on more complicated hardening laws we refer for instance to \cite{BroCarVal04, BroCarVal05, HofVal07}.
The elastoplastic initial-value constitutive model consists of the following components: \\

1. Additive decomposition of the strain tensor into the elastic and plastic parts:
\begin{equation}
\varepsilon=\varepsilon^e+\varepsilon^p.
\label{split}
\end{equation}

2. Linear elastic law between the stress and the elastic strain:
\begin{equation}
\sigma=\mathbb C\varepsilon^e,
\label{elastic_law}
\end{equation}
\hspace{1cm} where the fourth order tensor $\mathbb C$ is defined by (\ref{Hook_law}).

3. The von Mises yield function coupled with an isotropic hardening variable $\kappa$:
\begin{equation}
\Phi(\sigma,\kappa)=\sqrt{\frac{3}{2}}\|\mbox{dev}(\sigma)\|_F-(\sigma_y+H_m\kappa)\leq0,
\label{von_mises}
\end{equation}
\hspace{1cm} where $\sigma_y, H_m>0$ denote the initial yield stress and the hardening modulus, respectively.

4. The associated plastic flow rule:
\begin{equation}
\dot\varepsilon^p=\dot\gamma\frac{\partial \Phi}{\partial\sigma}=\dot\gamma\sqrt{\frac{3}{2}}\frac{\mbox{dev}(\sigma)}{\|\mbox{dev}(\sigma)\|_F},\quad \dot\gamma\geq0,
\label{flow_rule}
\end{equation}
\hspace{1cm} where $\dot\varepsilon^p$ and $\dot\gamma$ denote the time derivative of the plastic strain and the plastic multiplier, respectively.

5. The hardening law based on the accumulated plastic strain rate:
\begin{equation}
\dot\kappa=\sqrt{\frac{2}{3}}\|\dot\varepsilon^p\|_F=\dot\gamma.
\label{hardening_law}
\end{equation}
\hspace{1cm} Notice that the second equality in (\ref{hardening_law}) follows from (\ref{flow_rule}).

6. The loading/unloading conditions:
\begin{equation}
\dot\gamma\geq0,\;\; \Phi(\sigma,\kappa)\leq0,\;\; \dot\gamma\Phi(\sigma,\kappa)=0.
\label{load-unload}
\end{equation}

7. The initial conditions:
\begin{equation}
\varepsilon(x,t_0) = \varepsilon^e(x,t_0) = \varepsilon^p(x,t_0) = \sigma(x,t_0)=0, \;\;\kappa(x,t_0)=0,\quad x\in\Omega.
\end{equation}

The weak formulation of the corresponding elastoplastic problem can be found in \cite{HaRe99}. Here we will only consider a time discretized elastoplastic model.

\subsection{Time discretized elastoplastic model}

First, we derive an explicit form of the time discretized constitutive elastoplastic model. Then, we formulate the whole time discretized elastoplastic problem, similarly as in Subsection \ref{subsection_elast_model}.

Let us consider the following discretization of the time interval
$$t_0<t_1<\ldots<t_k<\ldots<t_N=t^*.$$
Let us denote $\sigma_k=\sigma_k(x)=\sigma(x,t_k)$, $x\in\Omega$ and similarly for other variables.
To approximate the time derivatives, we use the implicit Euler method. This method is often used in mathematical and engineering literature, see e.g. \cite{HaRe99,NPO08}. Other approximation schemes such as Crank-Nicholson scheme or discontinuous Galerkin are discussed e.g. in \cite{AC02}.
Then by (\ref{split}) and (\ref{elastic_law}),
\begin{equation}
\dot\varepsilon^p(t_{k+1})\approx \frac{\varepsilon^p_{k+1}-\varepsilon^p_k}{\triangle t_{k+1}}=\frac{\mathbb C^{-1}(\sigma_{k+1}^t-\sigma_{k+1})}{\triangle t_{k+1}},\quad \triangle t_{k+1}=t_{k+1}-t_{k},
\label{impl_euler}
\end{equation}
where
\begin{equation}
\sigma_{k+1}^t:=\sigma_{k}+\mathbb C\triangle\varepsilon_{k+1},\quad \triangle\varepsilon_{k+1}=\varepsilon_{k+1}-\varepsilon_{k}.
\label{trial}
\end{equation}
The stress tensor $\sigma_{k+1}^t$ is denoted as a trial stress tensor. By (\ref{split})--(\ref{trial}), we can formulate the time discretized elastoplastic {\it{constitutive}} problem as follows.
Given the values $\sigma_k$, $\kappa_k$, $\varepsilon_k$ of the stress, the isotropic hardening and the strain, respectively, at the time $t_k$ and given the incremental strain $\triangle\varepsilon_{k+1}$ for the interval $[t_k, t_{k+1}]$, solve the following system of algebraic equations
\begin{eqnarray}
\mathbb C^{-1}(\sigma_{k+1}^t-\sigma_{k+1})&=&\triangle\gamma_{k+1}\sqrt{\frac{3}{2}}\frac{\mbox{dev}(\sigma_{k+1})}{\|\mbox{dev}(\sigma_{k+1})\|_F}\label{flow_discr}\\
\kappa_{k+1}-\kappa_k&=&\triangle\gamma_{k+1}\label{hard_law_discr}
\end{eqnarray}
for the unknowns $\sigma_{k+1}$, $\kappa_{k+1}$, and $\triangle\gamma_{k+1}$, subject to the constraints
\begin{equation}
\triangle\gamma_{k+1}\geq0,\;\; \Phi(\sigma_{k+1},\kappa_{k+1})\leq0,\;\; \triangle\gamma_{k+1}\Phi(\sigma_{k+1},\kappa_{k+1})=0.
\label{load-unload_discr}
\end{equation}
This constitutive problem can be solved explicitly by the return mapping concept (see e.g. \cite{BLA99, NPO08}). It means that we firstly apply the elastic predictor, i.e., we verify whether $\Phi(\sigma_{k+1}^t,\kappa_{k})\leq0$. If it holds then
\begin{equation}
\triangle\gamma_{k+1}=0,\quad \sigma_{k+1}=\sigma_{k+1}^t, \quad \triangle \sigma_{k+1}=\mathbb{C}\triangle \varepsilon_{k+1},
\label{elast_pred}
\end{equation}
i.e., the stress increment satisfies the elastic law defined by \eqref{Hook_law}.
If $\Phi(\sigma_{k+1}^t,\kappa_{k})>0$, then by the plastic corrector we have
\begin{equation}
\triangle\gamma_{k+1}=\frac{1}{3\mu+H_m}\Phi(\sigma_{k+1}^t,\kappa_{k}),\quad \sigma_{k+1}=\sigma_{k+1}^t-\frac{3\mu}{3\mu+H_m}\sqrt{\frac{2}{3}} \Phi(\sigma_{k+1}^t,\kappa_{k}) \hat n(\sigma_{k+1}^t),
\label{plast_corr}
\end{equation}
where
\begin{equation}
\hat n(\tau)=\frac{\mbox{dev}(\tau)}{\|\mbox{dev}(\tau)\|_F}, \quad \tau \in S.
\label{hat_n}
\end{equation}
Notice that the second formula in (\ref{plast_corr}) is correctly defined since the denominator $\|\mbox{dev}(\sigma_{k+1}^t)\|_F>0$ for $\Phi(\sigma_{k+1}^t,\kappa_{k})>0$.
Let us define the stress and hardening operators $T_\sigma(\tau,\omega;.): S\rightarrow S$, $T_\kappa(\tau,\omega;.): S\rightarrow S$ with respect to parameters $\tau\in S$, $\omega\in\mathbb R_+$, such that for $\eta\in S$
\begin{eqnarray}
T_\sigma(\tau,\omega;\eta)&:=&
\mathbb C\eta-\frac{3\mu}{3\mu+H_m}\sqrt{\frac{2}{3}}\Phi^+(\tau+\mathbb C\eta,\omega)\hat n(\tau+\mathbb C\eta), \label{stress_op}\\
T_\kappa(\tau,\omega;\eta)&:=&
\frac{1}{3\mu+H_m}\Phi^+(\tau+\mathbb C\eta,\omega),\label{hardening_op}
\end{eqnarray}
respectively, where $\Phi^+$ denotes the positive part of the function $\Phi$. Then by (\ref{trial}), (\ref{hard_law_discr}), (\ref{elast_pred}), (\ref{plast_corr}), (\ref{stress_op}) and (\ref{hardening_op}),
\begin{equation}
\triangle\kappa_{k+1}=T_\kappa(\sigma_k,\kappa_k;\triangle\varepsilon_{k+1}),\quad \triangle\sigma_{k+1}=T_\sigma(\sigma_k,\kappa_k;\triangle\varepsilon_{k+1}).
\end{equation}
For the sake of brevity, we will denote the stress operator $T_\sigma(\sigma_k,\kappa_k;.)$ with respect to the current parameters $\sigma_k$ and $\kappa_k$ by $T_k(.)$. By \cite{BLA99, GV09, Sy09}, the operator $T_k: S\rightarrow S$ is potential, Lipschitz continuous, strongly monotone, and strongly semismooth on $S$.

Let us note that semismoothness was originally introduced by Mifflin \cite{Mi77} for functionals. Qi and J. Sun \cite{QS93} extended the definition of semismoothness to vector-valued function to investigate the superlinear convergence of the Newton method. The strong semismoothness of the Lipschitz continuous function $T_k(.)$ means that $T_k(.)$ is directionally differentiable on $S$ and has a quadratic approximate property at any $\eta\in S$, i.e.,
for any $\xi\in S$, $\xi\rightarrow0$, and any $T_k^o(\eta+\xi)\in\partial T_k(\eta+\xi)$,
\begin{equation}
T_k(\eta+\xi)-T_k(\eta)-T_k^o(\eta+\xi)\xi=O(\|\xi\|_F^2).
\label{semi}
\end{equation}
Here $\partial T_k(\eta+\xi)$ denotes the set of the Clark generalized derivatives of $T_k$ at $\eta+\xi$. Here we will choose the Clark generalized derivative $T_k^o$ of $T_k$ in the following way:

1. If $\Phi(\sigma_{k}+\mathbb C\eta,\kappa_{k})\leq0$, then
\begin{equation}
T_k^o(\eta)=\mathbb C.
\label{stress_op_deriv_elast}
\end{equation}

2. If $\Phi(\sigma_{k}+\mathbb C\eta,\kappa_{k})>0$, then
\begin{eqnarray}
T_k^o(\eta)&=&\mathbb C-\frac{3\mu}{3\mu+H_m}\sqrt{\frac{2}{3}}\hat n(\sigma_{k}+\mathbb C\eta)\otimes\frac{\partial\Phi(\sigma_{k}+\mathbb C\eta,\kappa_{k})}{\partial \eta}-\nonumber\\
&&-\frac{3\mu}{3\mu+H_m}\sqrt{\frac{2}{3}}\Phi(\sigma_{k}+\mathbb C\eta,\kappa_{k})\frac{\partial\hat n(\sigma_{k}+\mathbb C\eta)}{\partial \eta},
\label{stress_op_deriv_plast}
\end{eqnarray}
where
\begin{eqnarray*}
\frac{\partial\Phi(\sigma_{k}+\mathbb C\eta,\kappa_{k})}{\partial \eta}&=&2\mu\sqrt{\frac{3}{2}}\hat n(\sigma_{k}+\mathbb C\eta),\\
\frac{\partial\hat n(\sigma_{k}+\mathbb C\eta)}{\partial \eta}&=&2\mu\frac{\mathbb I_d-\hat n(\sigma_{k}+\mathbb C\eta)\otimes\hat n(\sigma_{k}+\mathbb C\eta)}{\|\mbox{dev}(\sigma_{k}+\mathbb C\eta)\|_F},\\
\mathbb I_d\xi&:=&\mbox{dev}(\xi),\quad\forall\xi\in S.
\end{eqnarray*}
Notice that $T_k$ is not differentiable at $\eta\in S$, $\Phi(\sigma_{k}+\mathbb C\eta,\kappa_{k})=0$. Otherwise $T_k^o(\eta)=\partial T_k(\eta)/\partial\eta$. By using (\ref{Hook_law})--(\ref{vol-dev_prop}),  (\ref{stress_op_deriv_elast}) and  (\ref{stress_op_deriv_plast}) we can derive the following uniform estimate for $T_k^o$  (see e.g. \cite{BLA99}):
\begin{equation}
\langle \mathbb C\xi,\xi\rangle\geq \langle T_k^o(\eta)\xi,\xi\rangle\geq \frac{H_m}{3\mu+H_m}\langle \mathbb C\xi,\xi\rangle\quad\forall \eta,\xi\in S,\;\forall k-\mbox{integer}.
\label{T_k^o-est}
\end{equation}

Let us recall that the stress, strain, hardening and displacement variables also depend on a spatial variable $x\in \Omega$. We consider the dependence of $T_k(\triangle\varepsilon_k)$ on $x$ in the following sense:
\begin{equation}
T_k(\triangle\varepsilon_k)=T_k(\triangle\varepsilon_k)(x):=T_\sigma(\sigma_k(x),\kappa_k(x);\triangle\varepsilon_k(x)).
\label{x-depend}
\end{equation}
Then we can substitute the stress operator $T_k$, defined by (\ref{stress_op}), into the balance equation (\ref{weak-balance}) to obtain the time discretized elastoplastic problem in the incremental form.

\begin{problem}[One time step elastoplastic problem in the incremental form]\label{problem2}

 {\it Given the stress field $\sigma_k\in [L^2(\Omega)]^{3\times 3}_{sym}$ and the isotropic hardening field $\kappa_k\in L^2(\Omega)$ at the time $t_k$, find the displacement $u_{k+1}=u_k+\triangle u_{k+1}\in V$, where the increment $\triangle u_{k+1}\in V$ solves the variational equation
\begin{equation}
\int_\Omega\langle T_k(\varepsilon(\triangle u_{k+1})),\varepsilon(v)\rangle_F dx=\int_\Omega \triangle g_{k+1}^Tv dx + \int_{\Gamma_N}\triangle F_{k+1}^Tv ds\quad\forall v\in V,
\label{weak-plast}
\end{equation}
with loading increments $\triangle F_{k+1}=F_{k+1}-F_k$, $\triangle g_{k+1}=g_{k+1}-g_k$.
Set the stress and isotropic hardening fields $\sigma_{k+1}=\sigma_{k}+\triangle\sigma_{k+1}$, $\kappa_{k+1}=\kappa_{k}+\triangle\kappa_{k+1}$ in the next time step $t_{k+1}$ from the relations
\begin{equation}
 \triangle\sigma_{k+1}=T_\sigma(\sigma_k,\kappa_k;\varepsilon(\triangle u_{k+1})),\quad
 \triangle\kappa_{k+1}=T_\kappa(\sigma_k,\kappa_k;\varepsilon(\triangle u_{k+1})),
\label{weak-plast2}
\end{equation}
almost everywhere in $\Omega$.
}
\end{problem}
Problem \ref{problem2} can be equivalently formulated as a minimization problem \cite{GV09, Sy09}. Since the operator $T_k$ is strongly monotone and Lipschitz continuous on $S$, the non-linear equation (\ref{weak-plast}) has a unique solution $\triangle u_{k+1}\in V$ (see e.g. \cite{FukKuf80}). As we will see in the next section, we will solve a linearized problem in each Newton iteration. To do this, it will be useful to define the bilinear form $a_{k}(u):V\times V\rightarrow \mathbb R$ for $u\in V$ by
\begin{equation}
a_{k}(u)(w,v)=\int_\Omega\langle T_k^o(\varepsilon(u))\varepsilon(w),\varepsilon(v)\rangle dx,\quad v,w\in V,
\end{equation}
where the operator $T_k^o(.)=T_k^o(.)(x)$ is defined by (\ref{stress_op_deriv_elast}) and (\ref{stress_op_deriv_plast}). 
The bilinear form $a_{k}$ is symmetric, bounded, and $V$-elliptic on $V$ due to (\ref{T_k^o-est}) and (\ref{a_e-est}).


\section{Semismooth Newton method in elastoplasticity} \label{section.Newton}

In this section, firstly, we approximate the time discretized elastoplastic problem by the finite element method and introduce the corresponding algebraic notation. Secondly, we introduce the semismooth Newton method for the problem.

\subsection{Finite element discretization and algebraic formulation}

Details to finite element implementation of elastoplastic problems can be found in \cite{CarKlo-2002} and \cite{BLA99, GV09}.

For the sake of simplicity, we assume a polyhedral 3D domain $\Omega$ and use the linear simplex elements. The corresponding shape regular triangulation is denoted by $\mathcal T_h$. Thus the space $V$ is approximated by its subspace $V_h$ of piecewise linear and continuous functions. Therefore the spaces of the strains, the stress and the isotropic hardening are approximated by piecewise constant functions.

Similarly to (\ref{weak-plast}), (\ref{weak-plast2}), we can formulate the one time step elastoplastic problem after the space discretization. Let $\sigma_{k,h}$, $\kappa_{k,h}$ be piecewise constant stress and hardening variables with respect to the triangulation $\mathcal T_h$ at the time $t_k$ obtained from a previous time process.
\begin{problem}[One time step elastoplastic problem in the incremental form after space discretization] {\it Find the displacement $u_{k+1,h}=u_{k,h}+\triangle u_{k+1,h}\in V_h$, where the increment $\triangle u_{k+1,h}\in V_h$ solves the variational equation
\begin{equation}
\int_\Omega\langle T_{k,h}(\varepsilon(\triangle u_{k+1,h})),\varepsilon(v_h)\rangle_F dx=\int_\Omega \triangle g_{k+1}^Tv_h dx + \int_{\Gamma_N}\triangle F_{k+1}^Tv_h ds\quad\forall v_h\in V_h,
\label{weak-plast_h}
\end{equation}
where $T_{k,h}(.):=T_\sigma(\sigma_{k,h},\kappa_{k,h};.)$. Set the stress and isotropic hardening fields $\sigma_{k+1,h}=\sigma_{k,h}+\triangle\sigma_{k+1,h}$, $\kappa_{k+1,h}=\kappa_{k,h}+\triangle\kappa_{k+1,h}$ in the next time step $t_{k+1}$ from the relations
\begin{equation}
\triangle\sigma_{k+1,h}=T_\sigma(\sigma_{k,h},\kappa_{k,h};\varepsilon(\triangle u_{k+1,h})),\quad
\triangle\kappa_{k+1,h}=T_\kappa(\sigma_{k,h},\kappa_{k,h};\varepsilon(\triangle u_{k+1,h}))
\label{weak-plast2_h}
\end{equation}
for every elements of $\mathcal{T}_h$.
}
\end{problem}
For the sake of simplicity, we do not consider finite element approximation of $F_k$ and $g_k$. Similarly as in (\ref{stress_op_deriv_elast}) and (\ref{stress_op_deriv_plast}), we can define the generalized derivative $T_{k,h}^o$ of $T_{k,h}$ and consequently also define the approximated bilinear form $a_{k,h}(u_h)$ for $u_h\in V_h$ by
\begin{equation}
a_{k,h}(u_h)(w_h,v_h)=\int_\Omega\langle T_{k,h}^o(\varepsilon(u_h))\varepsilon(w_h),\varepsilon(v_h)\rangle dx,\quad v_h,w_h\in V_h.
\label{bilinear_h}
\end{equation}

Each function $v_h=(v_{h,1},v_{h,2},v_{h,3})\in V_h$ can be represented by a vector $$\mbf v\in\mathbb R^n,\ \mbf v:=(v_{h,j}(x_i))_{i\in\{1,\ldots,\mathcal N\}, j\in\{1,2,3\}},$$ where $\mathcal N$ denotes the number of vertices of the triangulation $\mathcal T_h$ and $n=3\mathcal N$. The homogeneous Dirichlet boundary condition is represented by a restriction matrix $\mbf{B}_{U}\in\mathbb R^{m\times n}$, i.e.,
\begin{equation}
\mbf{B}_U\mbf{u}=\mbf o.
\end{equation}
Let $\mbf{R}_{T}\in\mathbb R^{12\times n}$ be a restriction operator for a displacement vector $\mbf u\in \mathbb R^n$ on a local element $T\in \mathcal T_h$, i.e.,
\begin{equation}
\mbf{u}_{T}=\mbf{R}_{T}\mbf u.
\label{displ_restr}
\end{equation}
We denote by $\mbf{u}_{k}$ and $\mbf{\triangle u}_{k+1}$ the displacement vector and the searching displacement increment at the time step $k$, respectively. We denote the load vector represented the volume and surface forces $F_k,g_k$ by $\mbf{f}_{k}$ and its increment $\mbf{\triangle f}_{k}$.

Further, we use a vector representation in $\mathbb R^6$ of the stress and strain tensors that is typical for an implementation of elastic problem, i.e.,
\begin{equation}
\mbf\sigma=(\sigma_{11},\sigma_{22},\sigma_{33},\sigma_{12},\sigma_{23},\sigma_{13})^T,\quad
\mbf\varepsilon=(\varepsilon_{11},\varepsilon_{22},\varepsilon_{33},2\varepsilon_{12},2\varepsilon_{23},2\varepsilon_{13})^{T}.
\label{stress_strain_vector}
\end{equation}
Notice that the stress and strain vectors have different structures in comparison to the above tensor notation. Therefore we must carefully distinguish this difference in algebraic representation of the operators $T_\sigma$, $T_\kappa$, $T_{k,h}$, and $T_{k,h}^o$. The vectors in sense of stress variables will be denoted by letters $\mbf\sigma$, $\mbf\tau$, the  vectors in sense of strain variables will be denoted by letters $\mbf\varepsilon$, $\mbf{\varepsilon}^{p}$, $\mbf\eta$, and $\mbf\xi$. Let $\mbf{\sigma}_{k,T}$ and $\mbf{\kappa}_{k,T}$ be the algebraic representation of $\sigma_{k,h}$ and $\kappa_{k,h}$ on an element $T\in\mathcal T_h$, respectively.

We introduce the algebraic representations  $\mbf C\in\mathbb R^{6\times 6}$, $\mbf E_{\varepsilon}\in\mathbb R^{6\times 6}$, $\mbf E_{\sigma}\in\mathbb R^{6\times 6}$, $\|\mbf.\|_{\sigma}$, $\mbf{\Phi}$, $\mbf{\hat n}$, $\mbf{T}_{\kappa,k,T}$, $\mbf{T}_{k,T}$, and $\mbf{T}^o_{k,T}$ of the Hooke tensor $\mathbb C$, the deviatoric operator $\mathbb I_d$ related to the strain and stress variables, the Frobenius norm with respect to a stress variable, the functions $\Phi$, $\hat n$, and the restrictions of the functions $T_{\kappa}\left( \sigma_k|_T, \kappa_k|_T, \cdot \right)$, $T_{k,h}$, $T^o_{k,h}$ on $T\in\mathcal T_h$, respectively, with respect to the vector form (\ref{stress_strain_vector}) of the stress and strain variables. The forms of matrices $\mbf C$, $\mbf{E}_{\varepsilon}$, $\mbf{E}_{\sigma}$, and the norm $\|\mbf.\|_{\sigma}$ are
$$
\begin{array}{rcl}
  \mbf C & := & \left[ \begin{array}{cccccc}
    \lambda+2\mu & \lambda      & \lambda      & 0   & 0   & 0 \\
    \lambda      & \lambda+2\mu & \lambda      & 0   & 0   & 0 \\
    \lambda      & \lambda      & \lambda+2\mu & 0   & 0   & 0 \\
    0            & 0            & 0            & \mu & 0   & 0 \\
    0            & 0            & 0            & 0   & \mu & 0 \\
    0            & 0            & 0            & 0   & 0   & \mu
  \end{array} \right], \\
  \mbf E_{\varepsilon} & := & \frac{1}{3} \left[ \begin{array}{rrrrrr}
     2 & -1 & -1 &   0 &   0 &   0\\
    -1 &  2 & -1 &   0 &   0 &   0\\
    -1 & -1 &  2 &   0 &   0 &   0\\
     0 &  0 &  0 & 1.5 &   0 &   0\\
     0 &  0 &  0 &   0 & 1.5 &   0\\
     0 &  0 &  0 &   0 &   0 & 1.5\\
  \end{array} \right], \\
  \mbf{E}_{\sigma} & := & \mbf{PE}_{\varepsilon},\quad \mbf P := \diag(1,1,1,2,2,2),
\end{array}
$$
and
$$
  \|\mbf{\tau}\|_{\sigma} := \left(\mbf{\tau}^T \mbf{P \tau}\right)^{1/2}, \quad \mbf{\tau} \in \mathbb{R}^6,
$$
respectively. Consequently the functions $\mbf{\Phi}$, $\mbf{\hat n}$, $\mbf{T}_{\kappa,k,T}$, $\mbf{T}_{k,T}$, $\mbf{T}^o_{k,T}$ are defined by \eqref{von_mises}, \eqref{hat_n}, \eqref{hardening_op}, \eqref{stress_op}, \eqref{stress_op_deriv_elast}, \eqref{stress_op_deriv_plast} in the following way:
$$
  \mbf{\Phi}(\mbf{\tau}, \mbf{\kappa}) := \sqrt{\frac{3}{2}} \| \mbf{E}_{\sigma} {\mbf{\tau}}\|_{\sigma} - \left(\sigma_y + H_m\mbf{\kappa}\right),
$$
$$
  \mbf{\hat n}(\mbf{\tau}) := \frac{\mbf{E}_{\sigma}\mbf{\tau}} {\|\mbf{E}_{\sigma}\mbf{\tau}\|_{\sigma}}, \quad \mbf{\overline{n}}_{k,T}(\mbf{\eta}) = \mbf{\hat n}\left( \mbf{\sigma}_{k,T} +\mbf{C\eta} \right),
$$
$$
  \mbf{T}_{\kappa,k,T}(\mbf{\eta}) := \frac{1}{3\mu + H_m}\mbf{\Phi}^{+}\left( \mbf{\sigma}_{k,T} +\mbf{C\eta}, \mbf{\kappa}_{k,T} \right),
$$
$$
  \mbf{T}_{k,T}(\mbf{\eta}) := \mbf{C\eta} - \frac{3\mu}{3\mu + H_m}\sqrt{\frac{2}{3}}\mbf{\Phi}^{+}\left( \mbf{\sigma}_{k,T} +\mbf{C\eta}, \mbf{\kappa}_{k,T} \right) \mbf{\overline n}_{k,T}(\mbf{\eta}),
$$
$$
  \mbf{T}^o_{k,T}(\mbf{\eta}) := \left\{\begin{array}{l}
    \mbf{C}, \quad \mbox{if } \mbf{\Phi}\left( \mbf{\sigma}_{k,T} +\mbf{C\eta}, \mbf{\kappa}_{k,T} \right)\leq 0, \mbox{ otherwise} \\
    \\
    \begin{array}{l}
      \mbf{C} - 2\mu\frac{3\mu}{3\mu + H_m} \mbf{E}_{\varepsilon} - \\
      - 2\mu\frac{3\mu}{3\mu + H_m} \sqrt{\frac{2}{3}} \frac{\sigma_y + H_m\mbf{\kappa}_{k,T}}{\|\mbf{E}_{\sigma}\left( \mbf{\sigma}_{k,T} +\mbf{C\eta} \right)\|_{\sigma}} \left( \mbf{\overline n}_{k,T}(\mbf{\eta}) \mbf{\overline n}^T_{k,T}(\mbf{\eta}) - \mbf{E}_{\varepsilon} \right),
    \end{array}
  \end{array}\right.
$$
respectively.

We also introduce the matrix $\mbf{G}_T\in\mathbb R^{6\times 12}$ representing the algebraical relation between the strain and the displacement (the exact form of $\mbf{G}_T$ is in \cite{AlbCarFun-2002}), i.e., the strain $\mbf{\varepsilon_T}$ on an element $T\in\mathcal T_h$ can be found by (\ref{displ_restr}) in the form
\begin{equation}
\mbf{\varepsilon}_{T}=\mbf{G}_T\mbf{R}_T\mbf{u}.
\label{stress_displ_alg}
\end{equation}

Based on the introduced notation we define the non-linear operator $\mbf{F}_{k}:\mathbb R^n\rightarrow\mathbb R^n$,
\begin{equation}
\mbf{F}_{k}(\mbf v)=\sum_{T\in\mathcal T_h}\left(\mbf{T}_{k,T}(\mbf{G}_T\mbf{R}_T\mbf{v})\right)^T\mbf{G}_T\mbf{R}_{T},\quad \mbf v\in \mathbb R^n,
\label{F_def}
\end{equation}
which plays the role of the left hand side in (\ref{weak-plast_h}).
Further we define the tangential and elastic stiffness matrices
\begin{eqnarray}
\mbf{K}_{k}(\mbf{v})&=&\sum_{T\in\mathcal T_h}\left(\mbf{T}^o_{k,T}({\mbf{G}_T}\mbf{R}_T\mbf{v}),{\mbf{G}_T}\mbf{R}_{T}\right)^T{\mbf{G}_T}\mbf{R}_{T},\quad \mbf{K}_{k}(\mbf{v})\in \mathbb R^{n\times n},
\; \mbf{v}\in\mathbb R^n,\label{K_def}\\
\mbf{K}_{e}&=&\sum_{T\in\mathcal T_h}\left(\mbf{C}{\mbf{G}_T}\mbf{R}_{T}\right)^T{\mbf{G}_T}\mbf{R}_{T},\quad \mbf{K}_{e}\in \mathbb R^{n\times n},
\end{eqnarray}
which represent the bilinear forms $a_{k,h}$ and $a_e$, respectively.
In particular, we denote the matrix $\mbf{K}_{k}(\mbf{\triangle u}_{k+1,i})$ briefly by $\mbf{K}_{k,i}$, where the reason of $\mbf{\triangle u}_{k+1,i}\in\mathbb R^n$ will be explained in the next section.

Let
$$
  \mbf{V} := \left\{\mbf{v}\in \mathbb{R}^n | \mbf{B}_U\mbf{v} = \mbf{o} \right\}.
$$
Then by using (\ref{F_def}), we can rewrite the equation (\ref{weak-plast_h}) as follows: find $\mbf{ \triangle u}_{k+1}\in\mbf{V}$ such that
\begin{equation}
\mbf{v}^T\left(\mbf{F}_{k}(\mbf{ \triangle u}_{k+1})-\mbf{\triangle f}_{k+1}\right)=0\quad\forall\mbf{v}\in\mbf{V},
\label{weak-plast_alg}
\end{equation}
where $\mbf{\triangle f}_{k+1}$ is the increment of the load vector.
Let $\mbf{\tilde u}_{k}\in\mathbb R^{n-m}$, $\mbf{\tilde f}_{k}\in\mathbb R^{n-m}$, $\mbf{\tilde K}_{k,i},\mbf{\tilde K}_{e}\in\mathbb R^{(n-m)\times(n-m)}$, and $\mbf{\tilde F}_{k}:\mathbb R^{n-m}\rightarrow\mathbb R^{n-m}$ denote the restrictions of $\mbf{u}_{k}$, $\mbf{ f}_{k}$, $\mbf{K}_{k,i},\mbf{K}_{e} $, and $\mbf{F}_{k}$ given by omitting the entries (degrees of freedom) corresponding to the prescribed Dirichlet boundary conditions. Then we can rewrite the equation (\ref{weak-plast_alg}) to the following system of non-linear equations:
\begin{equation}
\mbox{find } \mbf{\triangle u}_{k+1}\in\mbf{V}: \quad \mbf{\tilde F}_{k}(\mbf{ \triangle \tilde u}_{k+1})=\mbf{\triangle \tilde f}_{k+1}.
\label{non_system}
\end{equation}
The discretized elastoplastic problem can be solved by the following algorithm:

\begin{algorithm}[Solution of discretized elastoplastic problem]
\hspace{0.2cm}
\begin{spacing}{1.2}
\begin{algorithmic}[1]
  \STATE initial step: $\mbf{u}_{0} = \mbf o,\ \mbf{\varepsilon}_{0,T}=\mbf o,\ \mbf{\sigma}_{0,T} = \mbf o, \ \mbf{\kappa}_{0,T} = \mbf o$ for any $T\in\mathcal T_h$
  \FOR{$k = 0,\ \ldots, N-1$}
    \STATE find $\mbf{ \triangle u}_{k+1}\in\mbf{V}$: $\mbf{\tilde F}_k(\mbf{ \triangle \tilde u}_{k+1})=\mbf{\triangle \tilde f}_{k+1}$
    \FORALL{$T\in\mathcal T_h$}
      \STATE $\mbf{ \triangle \varepsilon}_{k+1,T}=\mbf{G}_T\mbf{R}_T \mbf{\triangle u}_{k+1},\;\mbf{ \varepsilon}_{k+1,T}=\mbf{\varepsilon}_{k,T}+\mbf{ \triangle \varepsilon}_{k+1,T}$
      \STATE $\mbf{ \triangle \sigma}_{k+1,T}=\mbf{T}_{k,T}(\mbf{ \triangle \varepsilon}_{k+1,T}),\;\mbf{ \sigma}_{k+1,T}=\mbf{\sigma}_{k,T}+\mbf{ \triangle \sigma}_{k+1,T}$
      \STATE $\mbf{ \triangle \kappa}_{k+1,T}=\mbf{T}_{\kappa,k,T}(\mbf{ \triangle \varepsilon}_{k+1,T}),\;\mbf{ \kappa}_{k+1,T} = \mbf{\kappa}_{k,T} +\mbf{ \triangle \kappa}_{k+1,T}$
    \ENDFOR
  \ENDFOR
\end{algorithmic}
\end{spacing}
\end{algorithm}

\subsection{Semismooth Newton method for one time step problem}

The non-linear system of equations (\ref{non_system}) is solved by the semismooth Newton method (see e.g. \cite{QS93}). The corresponding algorithm is following:

\begin{algorithm}[Semismooth Newton method]
\hspace{0.2cm}
\begin{spacing}{1.2}
\begin{algorithmic}[1]
  \STATE initialization: $\mbf{ \triangle u}_{k,0}=\mbf o$
  \FOR{$i=0,1,2,\ldots$}
    \STATE find $\mbf{\delta u}_{i}\in\mbf{V}$: $\mbf{\tilde K}_{k,i} \mbf{\delta \tilde u}_{i}=\mbf{\triangle \tilde f}_{k+1}-\mbf{\tilde F}_{k}(\mbf{ \triangle \tilde u}_{k,i})$
    \STATE compute $\mbf{ \triangle u}_{k,i+1}=\mbf{ \triangle u}_{k,i}+\mbf{\delta u}_{i}$
    \STATE {\bf if }{$\|\mbf{ \triangle u}_{k,i+1}-\mbf{ \triangle u}_{k,i}\|/(\|\mbf{ \triangle u}_{k,i+1}\|+\|\mbf{ \triangle u}_{k,i}\|)\leq\epsilon_{Newton}$} {\bf then stop}
  \ENDFOR
  \STATE set $\mbf{ \triangle u}_{k+1}=\mbf{ \triangle u}_{k,i+1}$
\end{algorithmic}
\end{spacing}
\end{algorithm}

Here $\epsilon_{Newton}>0$ is the relative stopping tolerance and $\mbf{\delta \tilde u_i} \in \mathbb{R}^{n-m}$ is the restriction of $\mbf{\delta u_i}$ given by omitting the entries (degrees of freedom) corresponding to the prescribed Dirichlet boundary conditions.
Since (\ref{T_k^o-est}) yields
\begin{equation}\label{spectral_equivalence}
\mbf{\tilde w}^T\mbf{\tilde K}_{e}\mbf{\tilde w}\geq \mbf{\tilde w}^T\mbf{\tilde K}_{k,i}\mbf{\tilde w}\geq\frac{H_m}{3\mu+H_m}\mbf{\tilde w}^T\mbf{\tilde K}_{e}\mbf{\tilde w}\quad\forall \mbf{\tilde w}\in\mbf{\mathbb{R}^{n-m}},\; \forall k,i-\mbox{integers}
\end{equation}
the matrices $\mbf{\tilde K}_{k,i}$ are spectrally equivalent to $\mbf{\tilde K}_{e}$. By \cite{KLV04}, it means that all types of preconditioners for elastic problems can be applied to the linearized problem in each Newton step as well.
If $H_m=0$, we obtain the elastic-perfectly plastic problem, and ${\tilde K}_{k,i}$ become generally singular. Such a problem is complicated not only for computations but also for the complete theory in terms of displacements. Therefore, the mentioned spectral equivalency seems to be problematic for too small values of $H_m$. However, most of  local stiffness matrices used in assembly of ${\tilde K}_{k,i}$ remain in an elastic mode and ${\tilde K}_{k,i}$ become more regular in real computations.
Indeed, we have successfully used the introduced algorithms even for perfectly plastic case \cite{CHS13}.


In \cite{BLA99, GV09, Sy09}, superlinear local convergence of the algorithm has been derived. Let us note that the convergence depends on the discretization parameter $h$ of the triangulation. Therefore we can expect that the finer the mesh, the bigger the number of the Newton iterations. In \cite{Sy09}, a damped semismoooth Newton method for such a problem has also been described. Such a method has again superlinear local convergence and additionally global convergence.


\section{TFETI method for solving linearized problems} \label{section.TFETI}
In this section, we describe the TFETI domain decomposition method for solving linearized problems appearing in Newton iterations and its optimal solver based on projected conjugate gradient method with preconditioning (PCGP). Finally, we summarize TFETI based algorithm for solving the whole elastoplastic problem.

\subsection{TFETI domain decomposition method}

In the previous section, we showed that we need to solve the following system of the linear equations:
\begin{equation}
  \mbox{find } \mbf{\delta u}_{i}\in\mbf{V}:\quad \mbf{\tilde K}_{k,i}\mbf{\delta \tilde u}_{i}=\mbf{\triangle \tilde f}_{k+1}-\mbf{\tilde F}_{k}(\mbf{ \triangle \tilde u}_{k,i})
\label{linear_system}
\end{equation}
in each time step $k$ and in each Newton iteration $i$. Since the TFETI domain decomposition method can be described independently of the indices $k,i$, we will schematically write the problem in the form:
\begin{equation}
  \mbox{find } \mbf{u}\in\mbf{V}:\quad \mbf{\tilde K\tilde u}=\mbf{\tilde f},
\label{linear_system2}
\end{equation}
where $\mbf{\tilde K}$, $\mbf{\tilde u}$, $\mbf{\tilde f}$ are the restriction of $\mbf{K,\ u,\ f}$ with respect to the Dirichlet boundary conditions respectively.
Let us note that $\mbf K$ is a symmetric and positive semidefinite matrix and $\mbf{\tilde K}$ is a symmetric and positive definite matrix. Therefore the linearized problem (\ref{linear_system2}) has a unique solution. The problem (\ref{linear_system2}) can be equivalently rewritten as a minimization problem:
\begin{equation}\label{linear_system3_minimum}
  \mbox{find } \mbf{u}\in\mbf{V}:\quad \mbf{J}(\mbf{u}) \leq \mbf{J}(\mbf{v}),\ \forall \mbf{v} \in \mbf{V},
\end{equation}
where
$$
  \mbf{J}(\mbf{v}) = \frac{1}{2}\mbf{v}^T\mbf{Kv} - \mbf{f}^T\mbf{v},\quad \mbf{v}\in\mbf{V}.
$$

The corresponding functional representation of (\ref{linear_system2}) is following: find $u_h\in V_h$ such that for any $v_h\in V_h$,
\begin{equation}
\int_\Omega\langle\mathbb D_h\varepsilon(u_h),\varepsilon(v_h)\rangle_F dx=\int_\Omega g^Tv_h dx + \int_{\Gamma_N}F^Tv_h ds-\int_\Omega\langle\tau_h,\varepsilon(v_h)\rangle_F dx,
\label{weak_discr}
\end{equation}
where $\mathbb D_h$ restricted on $T\in\mathcal T_h$ is a constant fourth order symmetric and elliptic tensor and $\tau_h$ restricted on ${T\in\mathcal T_h}$ is a constant second order tensor from $S$.
We can see that (\ref{weak_discr}) is related to (\ref{linear_system}), if we replace
$$\mathbb D_h\sim T^o_{k,h}(\triangle u_{k,h,i}),\;\;\tau_h\sim T_{k,h}(\triangle u_{k,h,i}),\;\;F\sim \triangle F_{k+1},\;\; g\sim \triangle g_{k+1},\;\;u_h\sim \delta u_{h,i}$$
and set $\triangle u_{k,h,i}$ and $\delta u_{h,i}$ as the functional representation of $\mbf{ \triangle u_{k,i}}$ and $\mbf{ \delta u_{i}}$, respectively.

Notice that the problem (\ref{weak_discr}) has a similar scheme as the elastic problem defined in Subsection \ref{subsection_elast_model}. Therefore the below introduced TFETI method can also be explained for elasticity in the same way.

The variational problem of the type (\ref{weak_discr}) can be equivalently formulated as a minimization problem:
\begin{equation}
\mbox{find } u_h\in V_h:\quad J_h(u_h)\leq J_h(v_h)\quad\forall v_h\in V_h,
\label{min_form}
\end{equation}
where
\begin{equation}
J_h(v_h)=\frac{1}{2}\int_\Omega\langle\mathbb D_h\varepsilon(v_h),\varepsilon(v_h)\rangle_F dx-\int_\Omega g^Tv_h dx - \int_{\Gamma_N}F^Tv_h ds+\int_\Omega\langle\tau_h,\varepsilon(v_h)\rangle_F dx.
\label{min_functional}
\end{equation}

\begin{figure}[n] 
\begin{center}
\includegraphics[height=5.5cm]{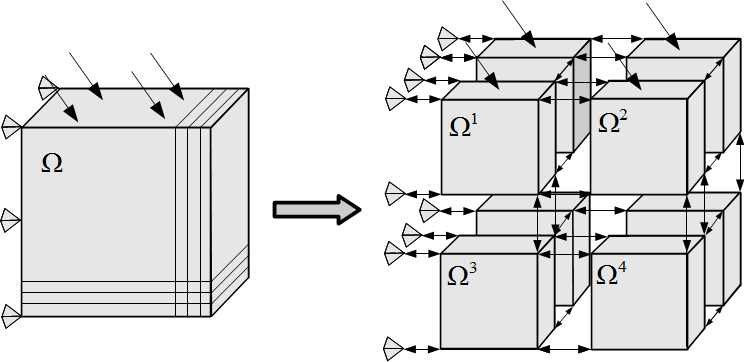}
\caption{\label{fig2}TFETI domain decomposition with subdomain renumbering}
\end{center}
\end{figure}

To apply the TFETI domain decomposition, we tear the body from the part of
the boundary with the Dirichlet boundary condition, decompose it into
subdomains, assign each subdomain by a unique number, and introduce new
``gluing'' conditions on the artificial intersubdomain boundaries and on the
boundaries with imposed Dirichlet condition (see Figure \ref{fig2}).

In particular, the polyhedral domain $\Omega$ is decomposed into a system of $s$ disjoint polyhedral subdomains $\Omega^p\subset \mathbb R^3$, $p=1,2,\ldots,s$.
We assume that the decomposition is consistent with the triangulation $\mathcal T_h$, i.e.,
\begin{equation}
\quad\forall T\in\mathcal T_h\; \exists !\ p\in\{1,2,\ldots,s\}:\quad T\subset\bar\Omega^p
\end{equation}
and define
\begin{equation}
\mathcal T_h^p:=\{T\in\mathcal T_h:\; T\subset\bar\Omega^p\},\quad \mathcal T_h=\bigcup_{p\in\{1,2,\ldots,s\}}\mathcal T_h^p.
\end{equation}
After the decomposition each boundary $\Gamma^p$ of $\Omega^p$ consists of three disjoint parts $\Gamma^p_U$, $\Gamma^p_N$, and $\Gamma^p_G$, $\Gamma^p = \overline\Gamma^p_U \cup \overline\Gamma^p_N \cup \overline\Gamma^p_G$, where
$$\Gamma^p_U=\Gamma_U\cap\Gamma^p,\;\;\Gamma^p_N=\Gamma_N\cap\Gamma^p,\;\;\Gamma^p_G=\bigcup_{q\in\{1,2,\ldots,s\}\setminus\{p\}}\Gamma^{pq}_G,$$
with $\Gamma^{pq}_G$ being the part of $\Gamma^p$ which is glued to $\Omega^q$, $p\neq q$.

Similarly to the definition of $V_h$, we can define the spaces $V_{h}^p$, $p=1,2,\ldots s$, of piecewise linear and continuous approximations of $H^1(\Omega^p)$:
\begin{equation} \label{def_V_h_p}
V_{h}^p:=\{v_h^p\in H^1(\Omega^p):\; v_h^p|_T\in [P_1]^3\;\;\forall T\in \mathcal T_h^p,\;\;v_h^p|_{\Gamma^p_U}=0\}.
\end{equation}
Let ${\tt V}_h:=V_{h}^1\times V_{h}^2\times\ldots\times V_{h}^s$ be a product space and
\begin{equation} \label{def_K_h}
{\tt K}_h:=\{{\tt v}_h=(v_h^1,\ldots,v_h^s)\in {\tt V}_h:\; v_h^p=v_h^q\;\;\mbox{on } \Gamma^{pq}_G\;\;\forall p,q\in\{1,2,\ldots,s\},\;p\neq q\}.
\end{equation}
Let us note that we slightly distinguish the notation introduced for the TFETI method from the notation introduced in Sections \ref{section.elastoplastic} and \ref{section.Newton}. For example, we write ${\tt V}_h$ for the TFETI method while in Section \ref{section.elastoplastic}, we used the notation $V_{h}$. A similar distinction is also introduced for the below algebraic description of the TFETI method.

Let
\begin{eqnarray}
{\tt J}_h({\tt v}_h)&=&\sum_{p=1}^s\left\{\frac{1}{2}\int_{\Omega^p}\langle\mathbb D_h\varepsilon({ v}_h^p),\varepsilon({ v}_h^p)\rangle_F dx-\int_{\Omega^p} g^T{ v}_h^p dx\right.\nonumber\\
&& \left.- \int_{\Gamma_N^p}F^T{ v}_h^p ds+\int_{\Omega^p}\langle\tau_h,\varepsilon({ v}_h^p)\rangle_F dx\right\}
\label{min_functional2}
\end{eqnarray}
be a functional defined on ${\tt V}_h$. Then the minimization problem (\ref{min_form}) can be equivalently rewritten into the form:
\begin{equation}
\mbox{find } {\tt u}_h\in {\tt K}_h:\quad {\tt J}_h({\tt u}_h)\leq {\tt J}_h({\tt v}_h)\quad\forall {\tt v}_h\in {\tt K}_h,
\label{min_form2}
\end{equation}
where ${\tt u}_h=\left((u_h)|_{\Omega^1},\ldots,(u_h)|_{\Omega^s}\right)$ and $u_h\in V_h$ solves (\ref{min_form}).

Each function ${\tt v}_h = \left(v_h^1,v_h^2, \ldots, v_h^s \right)$, ${\tt v}_h \in {\tt V}_h$, can be represented by a vector $\mathbf{v} \in \mathbb{R}^{\tt n}$, $\mathbf{v} = \left(\mbf{v}_1^T,\mbf{v}_2^T, \ldots, \mbf{v}_s^T \right)^T$, where $\mbf{v}_p \in \mathbb{R}^{n_p}, p\in \{1,2,\ldots, s\}$, is the algebraic representation of $v_h^p$ and ${\tt n} = \sum_{p=1}^s n_p$. Similarly we can find the vector $\mathbf{f} \in \mathbb{R}^{\tt n},\ \mathbf{f} = \left( \mbf{f}_1^T, \mbf{f}_2^T,\ldots, \mbf{f}_s^T \right)^T$, $\mbf{f}_p \in \mathbb{R}^{n_p}, p \in \{1,2,\ldots, s\}$, such that $\mbf{f}_p$ is the algebraic representation of the load restricted on $\Omega^p$ and $\Gamma_N^p$. Let the matrix $\mathbf{B}_G \in \mathbb{R}^{{\tt m}_G \times {\tt n}}$ represent the gluing conditions introduced in (\ref{def_K_h}) and  $\mathbf{B}_U \in \mathbb{R}^{{\tt m}_U \times {\tt n}}$ the Dirichlet boundary conditions introduced in (\ref{def_V_h_p}). Both matrices can be combined into one constraint matrix
\begin{equation}\label{equ.set.matrix.B}
  \mathbf{B} = \left[ \begin{array}{c}
    \mathbf{B}_G \\
    \mathbf{B}_U
  \end{array} \right], \quad \mathbf{B} \in \mathbb{R}^{{\tt m}\times {\tt n}},\quad {\tt m} = {\tt m}_G + {\tt m}_U.
\end{equation}
Typically ${\tt m}$ is much smaller than ${\tt n}$. Let us note that $\mathbf{B}$ can be assembled to have different forms: redundant, non-redundant or orthonormal. The rows of $\mathbf{B}_G$ in the standard redundant form are vectors of the order ${\tt n}$ with zero entries except 1 and -1 at appropriate positions. The orthonormal form is obviously obtained from the redundant one by only special treating of the rows of $\mathbf{B}$ corresponding to the nodes shared by more subdomains, i.e., for each of such nodes we take all corresponding rows in the redundant form, othonormalize them, and remove dependent ones. This can be done in parallel. For more details see~\cite{JPF97,PJF97,RF99,KW01a,FraPap-2003}. In fact all forms are applicable but due to simplicity of our presentation we use the orthonormal form of $\mathbf{B}$.

Let the matrix $\mathbf{K} \in \mathbb{R}^{{\tt n}\times {\tt n}}$, $\mathbf{K} = \diag(\mbf{K}_1, \mbf{K}_2, \ldots, \mbf{K}_s)$ denotes a symmetric positive semidefinite block diagonal matrix, where
$$
  \mbf{K}_p = \sum_{T\in\mathcal{T}_h^p} |T|\left(\mbf{D}_T \mbf{G}_T\mbf{R}_T^p \right)^T \mbf{G}_T\mbf{R}_T^p, \quad \mbf{K}_p \in \mathbb{R}^{n_p\times n_p}.
$$
Here $\mbf{D}_T \in \mathbb{R}^{6\times 6}$ is the algebraic representation of $\mathbb{D}_h|_T$ and $\mbf{R}_T^p \in \mathbb{R}^{12\times n_p}$ is a restriction operator for a displacement vector $\mbf u_p\in \mathbb R^{n_p}$ to a local element $T\in \mathcal T^p_h$. The diagonal blocks $\mbf{K}_p$, $p\in \{1,2,\ldots,s\}$, which correspond to the subdomains
$\Omega^p$, are positive semidefinite sparse matrices with known
kernels, the rigid body modes. 

The algebraical formulation of (\ref{min_form2}) is following:
\begin{equation} \label{equ.3radky}
  \left\{ \begin{array}{rcl}
  \mbox{find } \mathbf{u}\in\mathbf{V} \, & : & \, \mathbf{J}(\mathbf{u}) \leq \mathbf{J}(\mathbf{v})\quad \forall \mathbf{v}\in\mathbf{V}, \\
  \mathbf{J}(\mathbf{v}) &:=& \frac{1}{2}\mathbf{v}^T\mathbf{Kv} - \mathbf{f}^T\mathbf{v}, \\
  \mathbf{V} &:=& \left\{ \mathbf{v}\in\mathbb{R}^{\tt n}:\ \mathbf{Bv} = \mathbf{o} \right\}.
  \end{array} \right.
\end{equation}

Even though (\ref{equ.3radky}) is a standard convex quadratic programming problem,
its formulation is not suitable for numerical solution. The reasons are that
$\bf K$ is typically ill-conditioned, singular, and very large.

The complications mentioned above may be essentially reduced  by applying the
duality theory of convex programming (see, e.g., Dost\'al
\cite{DosOQPA-2007}), where all the constraints are enforced by the Lagrange multipliers $\vlambda$. The Lagrangian associated with problem (\ref{equ.3radky}) is
\begin{equation}
L({\bf v},\vlambda) = \mathbf{J}(\mathbf{v}) + \vlambda^{T} {\bf B v}.
\label{4.1}
\end{equation}
It is well known \cite{DosOQPA-2007} that (\ref{equ.3radky}) is
equivalent to the saddle point problem:
\begin{equation} \label{equ.min.max}
  \mbox{find } (\mathbf{u}, \boldsymbol{\lambda}) \in \mathbb{R}^{\tt n}\times \mathbb{R}^{\tt m}:\quad L(\mathbf{u},\boldsymbol{\nu})\leq L(\mathbf{u}, \boldsymbol{\lambda}) \leq L(\mathbf{v},\boldsymbol{\lambda})\quad \forall (\mathbf{v}, \boldsymbol{\nu})\in \mathbb{R}^{\tt n}\times \mathbb{R}^{\tt m}
\end{equation}
in sense that $\mathbf{u}$ solves \eqref{equ.3radky} if and only if $(\mathbf{u}, \boldsymbol{\lambda})$ solves \eqref{equ.min.max}.

\subsection{Optimal solvers to equality constrained problems}\label{section:solvers}

The solution of \eqref{equ.min.max} leads to the equivalent problem to find $({\mathbf{u}},{\boldsymbol{\lambda}})\in \dabR^{\tt n}\times \dabR^{\tt m}$ satisfying:
\begin{equation}\label{equ:sps}
     \left(
    \begin{array}{cc}
    \mathbf{K} & \mathbf{B}^{T}\\
    \mathbf{B} & \mathbf{0}
    \end{array}
    \right)
    \left(
    \begin{array}{c}
    \mathbf{u}\\
    \boldsymbol{\lambda}
    \end{array}
    \right)
    =
    \left(
    \begin{array}{c}
    \mathbf{f}\\
    \mathbf{o}
    \end{array}
    \right).
\end{equation}
The system \eqref{equ:sps} is uniquely solvable
which is guaranteed by the following necessary and sufficient conditions~\cite{BenGolLie05}:
\begin{eqnarray}
 \label{equ:NS1}
  && \Ker \mathbf{B}^{T} = \mathbf{o}, \\
 \label{equ:NS2}
  && \Ker \mathbf{K} \cap \Ker \mathbf{B} = \mathbf{o}. \label{equ.optimal.KerK}
\end{eqnarray}
Notice that \eqref{equ:NS1} is the condition on the full row-rank of $\mathbf{B}$.
Let us mention that an orthonormal basis of $\Ker \mathbf{K}$ is known a-priori and that its vectors
are columns of $\mathbf{R}\in \dabR^{{\tt n}\times l}$, $l={\tt n}-rank(\mathbf{K})$.

The first equation in \eqref{equ:sps} is satisfied if
\begin{equation} \label{equ:orth1}
    \mathbf{f}-\mathbf{B}^{T} {\boldsymbol{\lambda}} \in \Image \mathbf{K}
\end{equation}
and
\begin{equation} \label{equ:rec}
    {\mathbf{u}} = \mathbf{K}^\dagger (\mathbf{f}-\mathbf{B}^{T} {\boldsymbol{\lambda}}) + \mathbf{R} {\boldsymbol{\alpha}}
\end{equation}
for an appropriate ${\boldsymbol{\alpha}}\in \dabR^l$ and arbitrary matrix  $\mathbf{K}^\dagger$ satisfying $\mathbf{K} \mathbf{K}^\dagger \mathbf{K} = \mathbf{K}$. Here $\mathbf{K}^\dagger$ is a generalized inverse matrix whose application on a vector can be efficiently implemented (see Remark \ref{rem_pseudo}).

\begin{remark}\label{rem_pseudo}
The diagonal block $\mbf{K}_p$, $p\in \{1,2,\ldots,s\}$, which corresponds to the subdomain
$\Omega^p$, is positive semidefinite sparse matrix with known kernel basis created by the rigid body modes.
The first note about practical implementation of the action of the generalized inverse based on modified Cholesky factorization can be found in \cite{fr91}. This approach does not solve how to detect carefully zero pivots. In contrast to \cite{fr91} we use the kernel basis and fixing nodes strategy to effectively
regularize all blocks without extra fill in and then decompose them using any standard sparse Cholesky type factorization method for nonsingular matrices \cite{KucKozMar-2013, BrzDosKovKozMar-2010}.
The action of $\mathbf{K}^\dagger$ on a vector is naturally parallelized with respect to the subdomains and computed using backward and forward substitutions.
\end{remark}
\begin{remark}\label{rem_pseudo_standa}
Notice that the pseudoinverse $\mathbf{K}^\dagger$ can also be applied on vectors which do not belong to $\Image(\mathbf{K})$. Therefore we can write $\mathbf{K}^\dagger (\mathbf{f}-\mathbf{B}^T\boldsymbol{\lambda})=\mathbf{K}^\dagger \mathbf{f} - \mathbf{K}^\dagger (\mathbf{B}^T\boldsymbol{\lambda})$ in the other text. The random error corresponding with this partition is neglicted with respect to the implementation of $\mathbf{K}^\dagger$ introduced in Remark \ref{rem_pseudo}.
\end{remark}
The condition \eqref{equ:orth1} can be equivalently
written as
\begin{equation} \label{equ:orth2}
    \mathbf{R}^{T} (\mathbf{f}-\mathbf{B}^{T} {\boldsymbol{\lambda}}) = \mathbf{o}.
\end{equation}
Further substituting \eqref{equ:rec} into the second equation in \eqref{equ:sps}
we arrive at
\begin{equation} \label{equ:subst}
    -\mathbf{BK}^\dagger \mathbf{B}^{T} {\boldsymbol{\lambda}} + \mathbf{B R} {\boldsymbol{\alpha}} = - \mathbf{B K}^\dagger \mathbf{f}.
\end{equation}
Summarizing \eqref{equ:subst} and \eqref{equ:orth2} we find that the pair
$({\boldsymbol{\lambda}},{\boldsymbol{\alpha}})\in \dabR^m \times \dabR^l$ satisfies:
\begin{equation}\label{equ:sps3}
   \left(
    \begin{array}{cc}
    \mathbf{F} & \mathbf{N}^{T} \\
    \mathbf{N} & \mathbf{0}
    \end{array}
    \right)
    \left(
    \begin{array}{c}
    \boldsymbol{\lambda}\\
    \boldsymbol{\alpha}
    \end{array}
    \right)
    =
    \left(
    \begin{array}{c}
    \mathbf{d}\\
    \mathbf{e}
    \end{array}
    \right),
\end{equation}
where $\mathbf{F}:=\mathbf{BK}^\dagger \mathbf{B}^{T}$, $\mathbf{N}:=-\mathbf{R}^{T} \mathbf{B}^{T}$, $\mathbf{d}:=\mathbf{BK}^\dagger \mathbf{f}$, and $\mathbf{e}:=-\mathbf{R}^{T} \mathbf{f}$

Since $\mathbf{N}$ is of full row-rank as follow from \eqref{equ.optimal.KerK}, the inverse $(\mathbf{NN}^{T})^{-1}$ exists and $\mathbf{P_N}:=\mathbf{I}-\mathbf{N}^{T}(\mathbf{NN}^{T})^{-1}\mathbf{N}$ is well defined and represents the orthogonal projector onto $\Ker \mathbf{N}$. Applying $\mathbf{P_N}$ on the first equation in \eqref{equ:sps3} and checking that
$\mathbf{P_N} \mathbf{N}^{T} \boldsymbol{\alpha}=\mathbf{o}$ we eliminate $\boldsymbol{\alpha}$ and obtain that ${\boldsymbol{\lambda}}$ satisfies:
\begin{equation} \label{equ:projsps3}
    \mathbf{P_N}\mathbf{F} \boldsymbol{\lambda} = \mathbf{P_N} \mathbf{d}, \ \ \ \mathbf{N}\boldsymbol{\lambda} = \mathbf{e}.
\end{equation}
In practical computations, we further decompose $ {\boldsymbol{\lambda}}={\boldsymbol{\lambda}}_{\Image}+{\boldsymbol{\lambda}}_{\Ker}$   into two orthogonal components ${\boldsymbol{\lambda}}_{\Image}\in \Image \mathbf{N}^{T}$ and
${\boldsymbol{\lambda}}_{\Ker}\in \Ker \mathbf{N}$, substitute them in to \eqref{equ:projsps3} and get the problem
\begin{equation} \label{equ.proj.equation}
\mathbf{P_NF}\boldsymbol{\lambda}_{\Ker}=\mathbf{P_N}\left(\mathbf{d} - \mathbf{F}\boldsymbol{\lambda}_{\Image} \right) \mbox{ on } \Ker \mathbf{N},
\end{equation}
with $\boldsymbol{\lambda}_{\Image} = \mathbf{N}^T\left(\mathbf{NN}^T \right)^{-1}\mathbf{e}$, where $\mathbf{NN}^T$ is sparse and can be treated efficiently by parallel direct solvers such as MUMPS \cite{MUMPS}. Equation \eqref{equ.proj.equation} is solved efficiently by the projected conjugate gradient method with preconditioning (PCGP). In numerical examples, we consider both the computationally less expensive lumped preconditioner to $\mathbf{F}$ in the form $\overline{\mathbf{F}^{-1}}:=\mathbf{BKB}^T$ and the (quasi-) optimal Dirichlet preconditioner in the form $\overline{\mathbf{F}^{-1}}:=\overline{\mathbf{B}}\,\overline{\mathbf{S}}\,\overline{\mathbf{B}}^T$ introduced in \cite{fmr94}, where $\overline{\mathbf{S}}$ is the Schur complement to the block of $\mathbf{K}$ corresponding to the interior nodes and $\overline{\mathbf{B}}$ is the restriction of $\mathbf{B}$ to the remaining nodes, i.e., to the nodes with imposed Dirichlet and gluing conditions. The condition number estimates contain a linear factor of $H/h$ ($H$ is the domain decomosition step and $h$ is the discretization step) for the lumped preconditioner as opposed to a polylogarithmic bound for the Dirichlet
preconditioner. More details about other forms of the Dirichlet preconditioners and estimates of condition number can be found in \cite{RF99,KW01a}.
We obtain the following algorithmic scheme for the solution of \eqref{equ.3radky}:

\begin{algorithm}[Linear solver based on the TFETI method] \label{algorithm_TFETI}
\hspace{0.2cm}
\begin{spacing}{1.2}
\begin{algorithmic}[1]
  \STATE Set $\mathbf{N}:=-\mathbf{R}^{T} \mathbf{B}^{T}$, $\mathbf{H}:=(\mathbf{NN}^{T})^{-1}$, $\mathbf{d}:=\mathbf{BK}^\dagger \mathbf{f}$, and $\mathbf{e}:=-\mathbf{R}^{T} \mathbf{f}$.
  \STATE Compute ${\boldsymbol{\lambda}}_{\Image} := \mathbf{N}^{T} \mathbf{H e}$.
  \STATE Set $\tilde{\mathbf{d}}:=\mathbf{d}-\mathbf{F}{\boldsymbol{\lambda}}_{\Image}$.
  \STATE Compute ${\boldsymbol{\lambda}}_{\Ker}$ from \eqref{equ.proj.equation} by PCGP:
  \STATE \hspace{0.4cm} $\mathbf{r}^{0}=\tilde{\mathbf{d}}$, $\ \boldsymbol{\lambda}_{\Ker}^{0}=\mathbf{o}.$
  \STATE \hspace{0.4cm} {\bf for }$j=1, 2, \ldots,$ until convergence {\bf do}
    \STATE \hspace{0.8cm} Project $\mathbf{w}^{j-1}=\mathbf{P_Nr}^{j-1}.$
    \STATE \hspace{0.8cm} Precondition $\mathbf{z}^{j-1}=\overline{\mathbf{F}^{-1}}\mathbf{w}^{j-1}.$
    \STATE \hspace{0.8cm} Re-project $\mathbf{y}^{j-1}=\mathbf{P_Nz}^{j-1}.$
    \STATE \hspace{0.8cm} $\boldsymbol{\beta}^{j}=(\mathbf{y}^{j-1})^T \mathbf{w}^{j-1}/(\mathbf{y}^{j-2})^T \mathbf{w}^{j-2} \qquad(\boldsymbol{\beta}^{1}=0).$
    \STATE \hspace{0.8cm} $\mathbf{p}^{j}=\mathbf{y}^{j-1}+\boldsymbol{\beta}^{j}\mathbf{p}^{j-1}\hspace{3.0cm}(\mathbf{p}^{1}=\mathbf{y}^{0}).$
    \STATE \hspace{0.8cm} $\boldsymbol{\gamma}^{j}=(\mathbf{y}^{j-1})^T \mathbf{w}^{j-1}/(\mathbf{p}^{j})^T \mathbf{Fp}^{j}.$
    \STATE \hspace{0.8cm} $\boldsymbol{\lambda}_{\Ker}^{j}=\boldsymbol{\lambda}_{\Ker}^{j-1}+\boldsymbol{\gamma}^{j}\mathbf{p}^{j}.$
    \STATE \hspace{0.8cm} $\mathbf{r}^{j}=\mathbf{r}^{j-1}-\boldsymbol{\gamma}^{j}\mathbf{Fp}^{j}.$
    \STATE \hspace{0.8cm} {\bf if }$\|\mathbf{w}^{j-1}\|\leq\epsilon_{PCGP} \|r^0\|$ {\bf then stop.}
  \STATE \hspace{0.4cm} {\bf end for}
  \STATE \hspace{0.4cm} $\boldsymbol{\lambda}_{\Ker}=\boldsymbol{\lambda}_{\Ker}^j$.
  \STATE Set ${\boldsymbol{\lambda}}:={\boldsymbol{\lambda}}_{\Image}+{\boldsymbol{\lambda}}_{\Ker}$.
  \STATE Compute ${\boldsymbol{\alpha}}:=\mathbf{HN}(\mathbf{d}-\mathbf{F}{\boldsymbol{\lambda}})$.
  \STATE Compute ${\mathbf{u}}:=\mathbf{K}^\dagger(\mathbf{f}-\mathbf{B}^{T}{\boldsymbol{\lambda}} )+\mathbf{R}{\boldsymbol{\alpha}}$.
\end{algorithmic}
\end{spacing}
\end{algorithm}

\begin{remark}
Action of $\mathbf{H}$ on a vector may be efficiently implemented by the sparse Cholesky factorization of $\mathbf{NN}^T$.
\end{remark}

\subsection{TFETI based algorithms for solving elastoplastic problem}

In this subsection, we summarize the above proposed algorithms for solving elastoplastic problems and modify them with respect to the use of the TFETI domain decomposition method.

\begin{algorithm}[Algorithm for solving elastoplastic problem - sequential version] \label{TFETI_elasto_plastic_algorithm}
\hspace{0.2cm}
\begin{spacing}{1.2}
\begin{algorithmic}[1]
  \STATE $\mathbf{u}_{0} = \mathbf{o}, \boldsymbol{\sigma}_{0,T} = \mathbf{o}, \boldsymbol{\kappa}_{0,T} = \mathbf{o},\ T\in\mathcal{T}_h$ (initial step)
  \FOR{$k = 0,\, 1,\, 2,\,  \ldots,\, N-1$  (time steps)}
    \STATE set $\triangle\mathbf{u}_{k+1,0} = \mathbf{o}$ (zero approximation)
    \FOR{$i=1,2,\ldots$  (Newton iterations)}
      \FOR{ $p = 1,2,\ldots,s $ (cycle over subdomains)}
        \STATE restrict $\triangle\mathbf{u}_{k+1,i-1}, \triangle \mathbf{f}_{k+1}$ into subdomain variables $\mbf{\triangle u}_{k+1,i-1}^{p}$, $\mbf{\triangle f}_{k+1}^{p}$
        \STATE call Algorithm \ref{algorithm.compute_variables} with ($\mbf{\triangle u}_{k+1,i-1}^{p},\mbf{\triangle f}_{k+1}^p,\, \boldsymbol{\sigma}_{k,T},\, \boldsymbol{\kappa}_{k,T}$, $T\in\mathcal{T}_h^p$) to find output variables $\mbf{K}^{p}_{k,i}$, $\mbf{f}^{p}_{k,i}$, $\triangle\boldsymbol{\sigma}_{k+1,T}$, $\triangle\boldsymbol{\kappa}_{k+1,T}$, $T\in\mathcal{T}_h^p$.
        \STATE collect $\mbf{K}^{p}_{k,i}\, \mbf{f}^{p}_{k,i},$ into global variables $\mathbf{K}_{k,i}$, $\mathbf{f}_{k,i}$
      \ENDFOR\ (cycle over subdomains)
      \STATE solve by Algorithm \ref{algorithm_TFETI} the problem: $\mbox{find } \boldsymbol{\delta}\mathbf{u}_i \in \mathbf{V}$ such that
      $$\mathbf{J}_{k,i}(\boldsymbol{\delta}\mathbf{u}_i) \leq \mathbf{J}_{k,i}(\mathbf{v})\ \ \forall \mathbf{v}\in \mathbf{V}, \mbox{ with } \mathbf{J}_{k,i}(\mathbf{v}) = \frac{1}{2}\mathbf{v}^T\mathbf{K}_{k,i}\mathbf{v} - \mathbf{f}_{k,i}^T\mathbf{v}$$
      \STATE $\triangle\mathbf{u}_{k+1,i} = \triangle\mathbf{u}_{k+1,i-1} + \mbf{\delta}\mathbf{u}_{i} $ (displacement update)
      \STATE {\bf if} $\|\triangle\mathbf{u}_{k+1,i} - \triangle\mathbf{u}_{k+1,i-1}\| / \left( \|\triangle\mathbf{u}_{k+1,i}\| + \| \triangle\mathbf{u}_{k+1,i-1}\| \right) \leq \epsilon_{Newton}$ {\bf then stop}
    \ENDFOR (Newton iter.)
    \STATE $\mathbf{u}_{k+1} = \mathbf{u}_{k} + \triangle\mathbf{u}_{k+1,i},$
    \STATE $\boldsymbol{\sigma}_{k+1,T} = \boldsymbol{\sigma}_{k,T} + \triangle\boldsymbol{\sigma}_{k+1,T},\ \boldsymbol{\kappa}_{k+1,T} = \boldsymbol{\kappa}_{k,T} + \triangle\boldsymbol{\kappa}_{k+1,T},\, T\in\mathcal{T}_h$
  \ENDFOR\ (cycle over time steps)
\end{algorithmic}
\end{spacing}
\end{algorithm}

The assembling procedure for subdomain data looks as follows.

\begin{algorithm}[Assemble all data corresponding to a subdomain $\Omega^{p}$] \label{algorithm.compute_variables}
\hspace{0.2cm}
\begin{spacing}{1.2}
\begin{algorithmic}[1]
  \STATE Input: $\mbf{\triangle u}_{k+1,i-1}^{p},\, \mbf{\triangle f}_{k+1}^{p},\, \boldsymbol{\sigma}_{k,T},\, \boldsymbol{\kappa}_{k,T},\, T\in\mathcal{T}_h^p$.
  \STATE $\mbf{f}_{k,i}^p = \mbf{\triangle f}_{k+1}^p$
  \STATE $\mbf{K}^p_{k,i} = \mbf{O}$
  \FOR{$T\in \mathcal{T}_h^p$}
    \STATE compute $|T|$ (volume of the element $T$)
    \STATE $\triangle\boldsymbol{\sigma}_{k+1,T} = \mbf{T}_{k,T}\left(\mbf{G}_T\mbf{R}_T^p\mbf{\triangle u}_{k+1,i-1}^p \right)$
    \STATE $\triangle\boldsymbol{\kappa}_{k+1,T} =  \mbf{T}_{\kappa,k,T}\left(\mbf{G}_T\mbf{R}_T^p\mbf{\triangle u}_{k+1,i-1}^p \right)$
    \STATE $\mbf{f}_{k,i}^p = \mbf{f}_{k,i}^p - |T|\left( \triangle\boldsymbol{\sigma}_{k,T} \right)^T \mbf{G}_T\mbf{R}_T^p$
    \STATE $\mbf{K}_{k,i}^p = \mbf{K}_{k,i}^p + |T| \left( \mbf{T}_{k,T}^o\left( \mbf{G}_T\mbf{R}_T^p\mbf{\triangle u}_{k+1,i-1}^p \right) \mbf{G}_T\mbf{R}_T^p \right)^T \mbf{G}_T\mbf{R}_T^p$
  \ENDFOR\ (cycle over elements)
  \STATE Output: $\mbf{K}^{p}_{k,i},\, \mbf{f}^{p}_{k,i},\, \triangle\boldsymbol{\sigma}_{k+1,T}, \triangle\boldsymbol{\kappa}_{k+1,T},\, T\in\mathcal{T}_h^p$
\end{algorithmic}
\end{spacing}
\end{algorithm}

Loop over subdomains and all subdomain operations may be implemented in parallel.
Parallelization of FETI/TFETI is based on distributing matrix portions among processing units. This allows algorithms to be almost the same in sequential and parallel versions; only data structure implementation differs. Most of computations (subdomain operations) appearing are purely local and therefore parallelizable without any data transfers.
Each of cores works with the local part associated with its subdomains. Natural effort using the massively parallel computers is to maximize the number of subdomains so that the sizes of the subdomain stiffness matrices are reduced. This accelerates their factorization and subsequent pseudoinverse application which belongs to the most time consuming action. On the other hand, negative effect of that is an increase of the null space dimension and the number of Lagrange multipliers. This leads to larger coarse problems, i.e., applications of the projector $\mathbf{P_N}$ which are scalable only up to a few thousands of cores and then the coarse problem solution starts to dominate. For the numerical solution of such large problems we recommend to use a hybrid FETI method \cite{KlaRhe-2010}.


\section{Numerical experiments}
The proposed algorithms were implemented in {\tt MatSol} library \cite{MatSol-2009} developed in Matlab and parallelized using Matlab Distributed Computing Server and Matlab Parallel Toolbox. For all computations we use maximum 28 cores with 2GB memory per core of the HP Blade system, model BLc7000. The numbers of subdomains are chosen to keep the number of nodes per subdomain approximately constant except for the coarsest mesh level.

\begin{figure}[n] 
\begin{minipage}[t]{0.45\textwidth}
  \center
   \includegraphics[width=\textwidth]{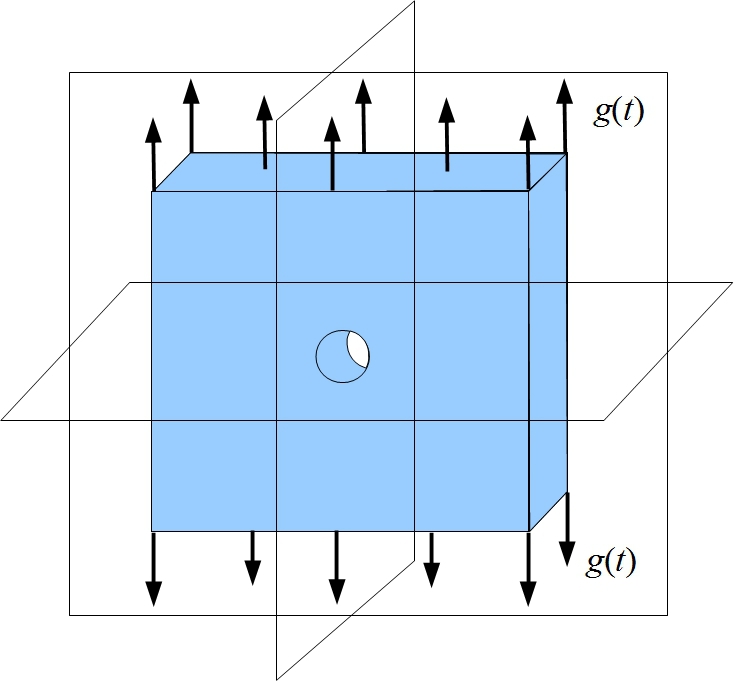}
   \caption{\small{Geometry of the whole body}}
   \label{fig04g}
\end{minipage}
\hfill
\begin{minipage}[t]{0.45\textwidth}
  \center
   \includegraphics[width=0.9\textwidth]{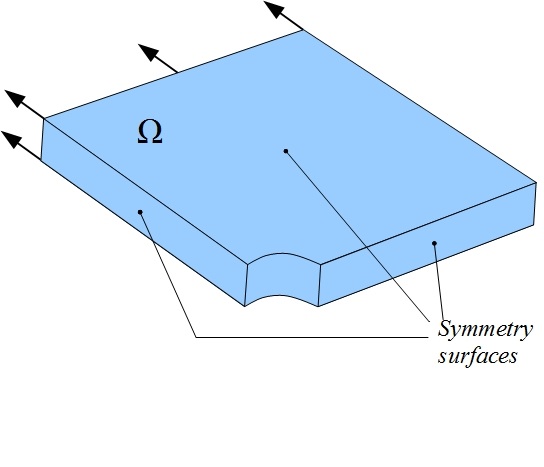}
   \caption{\small{Geometry of eighth of body}}
   \label{fig04g2}
\end{minipage}
\end{figure}

The performance is demonstrated on an elastoplastic homogeneous thin plate of sizes $20 \times 20 \times 2$ with the circular hole of radius 1 in the center. Its geometry with imposed boundary conditions and indicated symmetry planes is depicted in Figure \ref{fig04g}. Thus we consider only eighth of the thin plate in our computations (see Figure \ref{fig04g2}) and impose symmetry conditions on three symmetry surfaces. A similar benchmark was solved in \cite{Stein, GV09, CerKozMar-2011}.

\begin{figure}[htbp]
\begin{minipage}[t]{0.45\textwidth}
  \center
   \includegraphics[width=\textwidth]{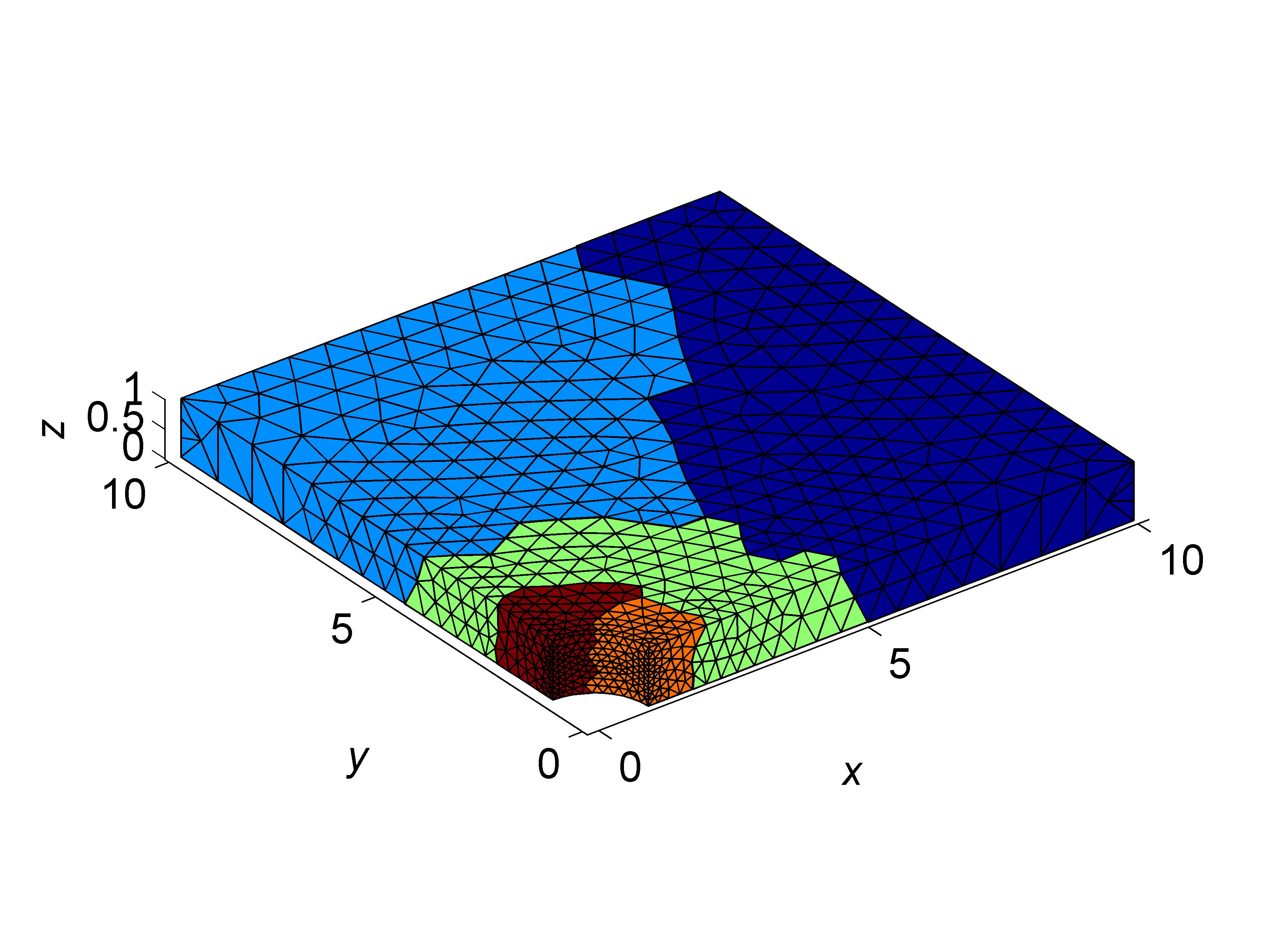}
   \caption{\small{Domain decomposition into 5 subdomains}}
   \label{fig04d}
\end{minipage}
\hfill
\begin{minipage}[t]{0.45\textwidth}
  \center
   \includegraphics[width=0.9\textwidth]{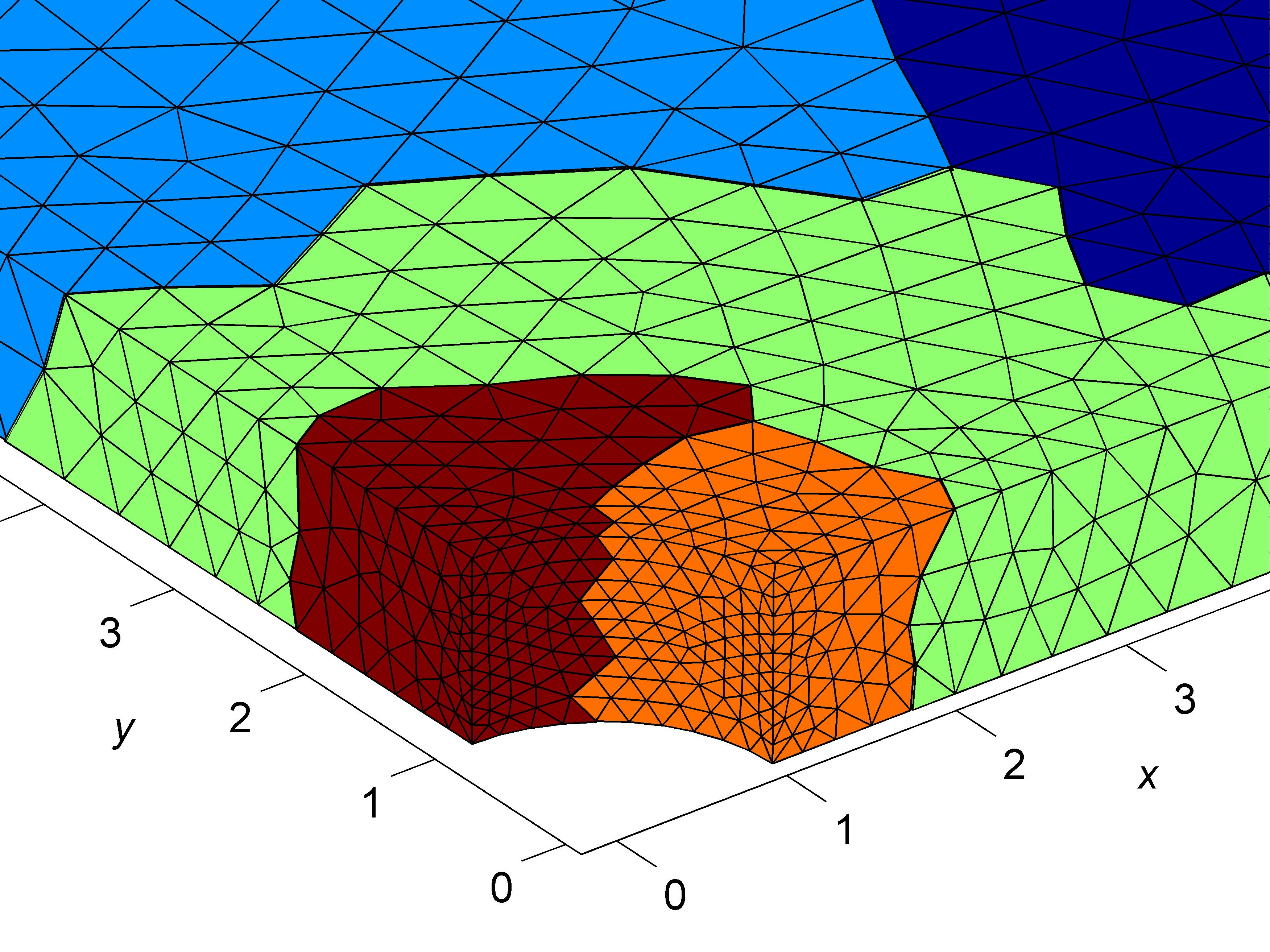}
   \caption{\small{Zoom of Figure \ref{fig04d} near the hole}}
   \label{fig04d2}
\end{minipage}
\end{figure}

The elastoplastic body $\Omega$ is made of homegeneous isotropic material with the parameters $$E = 206\,900, \, \nu = 0.29, \, \sigma_y = 450, \,\ \mbox{and}\ H_m = 10\,000,$$ where $E$ and $\nu$ are the Young modulus and the Poisson ratio, respectively, which are related with the Lam\'e coefficients $\lambda$ and $\mu$ by the following formulas: $$ \lambda = \frac{E\nu}{(1+\nu)(1-2\nu)}=110\,743.8,\quad \mu = \frac{E}{2(1+\nu)}=801\,938.$$
The indicated traction forces with the history of loading taking into account are prescribed by the function $$g(t) = 400\sin(2\pi t),\ t \in [t_0,t^*],\ t_0=0,\ t^*=\frac{1}{4}.$$
Since the elastoplastic model is rate-independent and any local unloading is not expected with respect to the prescribed load history, the results should be independent of the chosen time discretization.
Let us consider two variants of the equidistant time discretization characterized by the time step $\Delta t$:
\begin{itemize}
\item[$(a)$] $\Delta t=1/4$, $N=1$,
\item[$(b)$] $\Delta t=1/32$, $N=8$.
\end{itemize}
For the spatial discretization of $\Omega$, let us consider five levels of tetrahedral meshes generated by Ansys with
$$520, \, 1\, 623 , \, 9\, 947, \, 166\, 374, \, \ \mbox{and}\ 309\,546 $$ nodes decomposed into subdomains using Metis. An example of such decomposition is depicted in Figure \ref{fig04d} and a zoom near the hole in Figure \ref{fig04d2}.

In Table \ref{tab:diferent_meshes}, we report on numerical results for the time discretization $(a)$, different mesh levels, the Dirichlet and the lumped preconditioners.  The stopping tolerances of the Newton and the PCGP algorithms are
\begin{equation}\label{tolerances}
\epsilon_{Newton} = 10^{-4}\mbox{ and } \epsilon_{PCGP} = 10^{-7},
\end{equation}
respectively.
We see that the number of the Newton iterations remains almost constant for all meshes. Similar behavior was observed in the numerical results by \cite{GV09}. Note that even the total  number of the PCGP iterations (which means the sum of PCGP iterations over all Newton steps) increases only moderately for finer meshes and both preconditioners. We observe similarly as in ~\cite{fmr94} that the Dirichlet preconditioner is worse because of medium size of our problems, i.e., the total number of PCGP iterations is less than for the lumped preconditioner but not sufficiently.

In Table \ref{tab:martin5_upr}, the number of the PCGP iterations, the number of plastic elements, and the relative convergence criterion are reported for each Newton's iteration of Algorithm \ref{TFETI_elasto_plastic_algorithm}. In this case, the time discretization $(a)$, the finest mesh level, and the lumped preconditioner are considered together with sharper stopping tolerances
\begin{equation}\label{tolerances2}
\epsilon_{Newton} = 10^{-9}\mbox{ and } \epsilon_{PCGP} = 10^{-14},
\end{equation}
so that more Newton's iterations can be studied. Since the initial Newton iteration is taken as zero vector (of displacements), the first iteration of the Newton method actually solves a purely elastic problem. The second and further iterations already reflect an elastoplastic behaviour with stabilizing number of plastic elements and higher number of PCGP iterations. This is caused by larger condition number of a considered matrix in accordance with estimate \eqref{spectral_equivalence}. Once plastic elements are identified in the 7th step, we observe the quadratic convergence of the Newton method. Such behavior agrees with the theoretical results in \cite{BLA99}.

\begin{figure}[htbp]
\begin{minipage}{0.48\textwidth}
	\center \includegraphics[width=\textwidth]{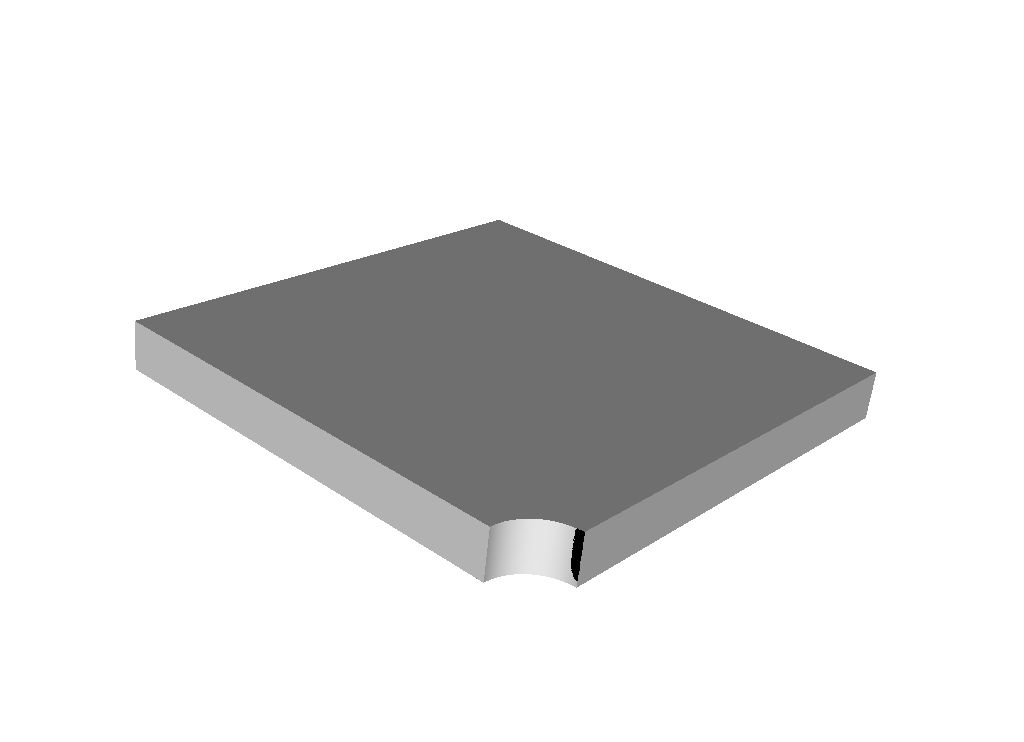}
	\caption{\small{Elastic (gray color) and plastic (black color) elements at time $t_2$}}
	\label{fig11}
\end{minipage}
\hfill
\begin{minipage}{0.48\textwidth}
  \center
   \includegraphics[width=\textwidth]{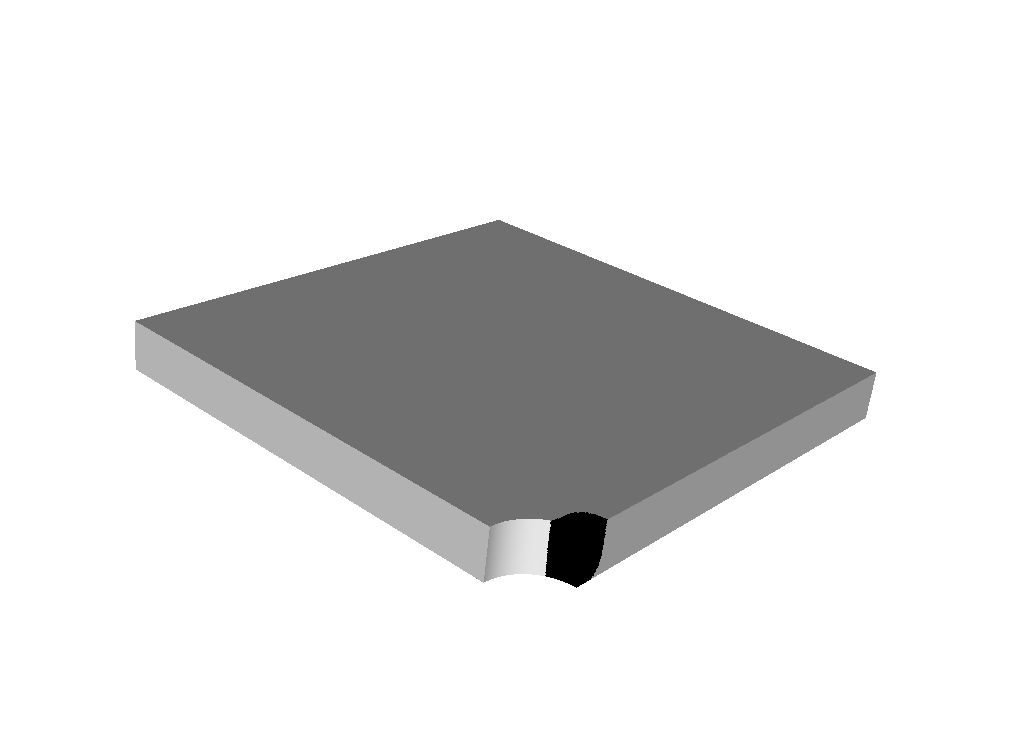}
   \caption{\small{Elastic (gray color) and plastic (black color) elements at time $t_4$}}
   \label{fig12}
\end{minipage}
\begin{minipage}{0.48\textwidth}
  \center
   \includegraphics[width=\textwidth]{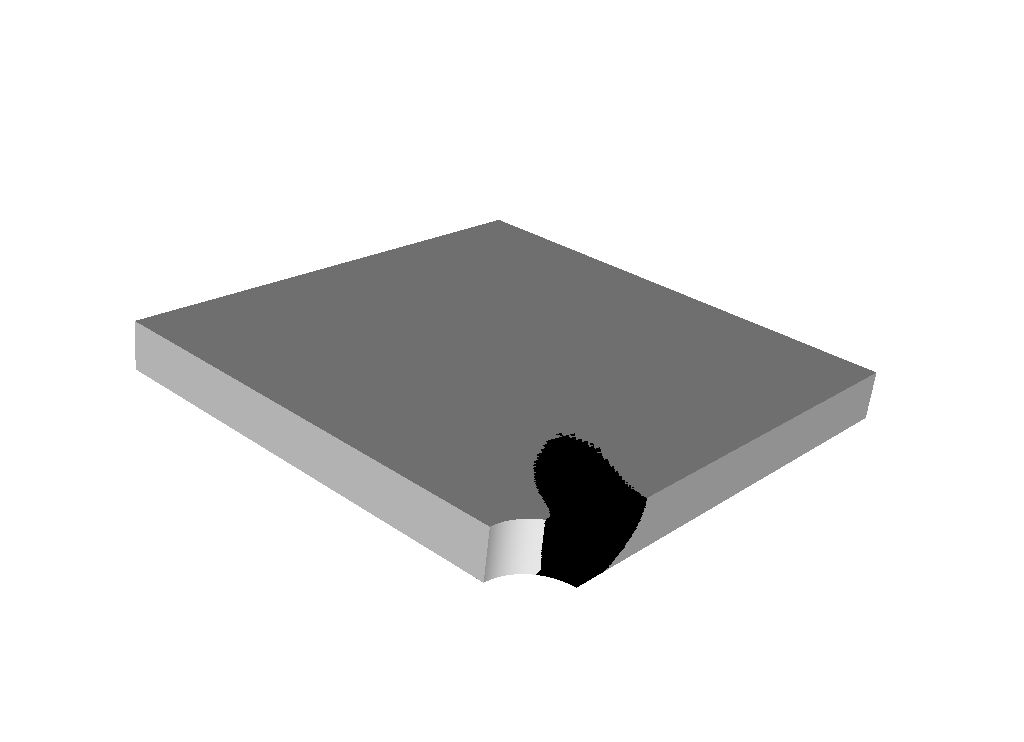}
   \caption{\small{Elastic (gray color) and plastic (black color) elements at time $t_6$}}
   \label{fig13}
\end{minipage}
\hfill
\begin{minipage}{0.48\textwidth}
  \center
   \includegraphics[width=\textwidth]{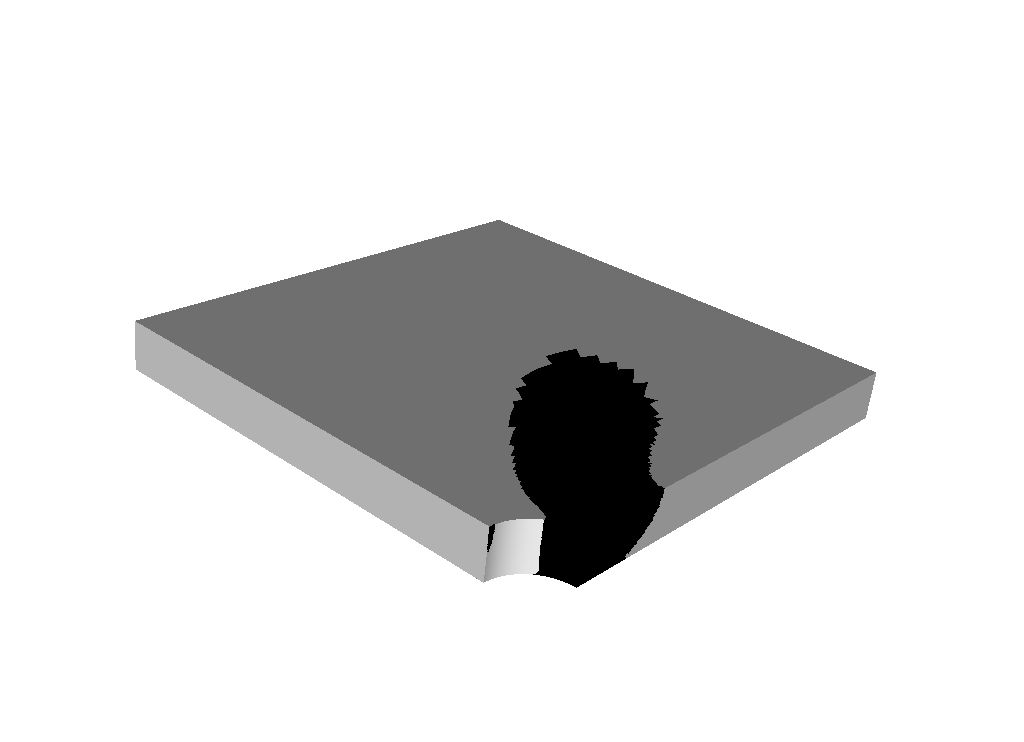}
   \caption{\small{Elastic (gray color) and plastic (black color) elements at time $t_8$}}
   \label{fig14}
\end{minipage}
\vspace{1cm} \\
\begin{minipage}{0.48\textwidth}
  \center
   \includegraphics[width=\textwidth]{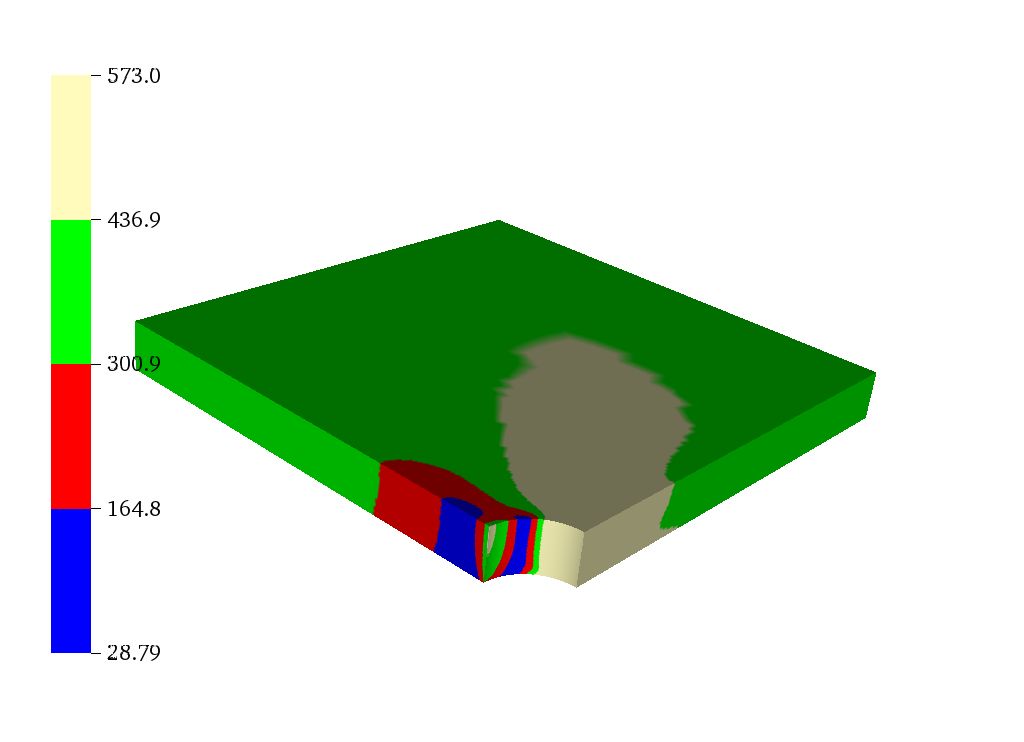}
   \caption{\small{Distribution of von Mises stress $\|\mbox{dev}(\sigma)\|_F$ at $t_8$}}
   \label{fig15}
\end{minipage}
\hfill
\begin{minipage}{0.48\textwidth}
  \center
   \includegraphics[width=\textwidth]{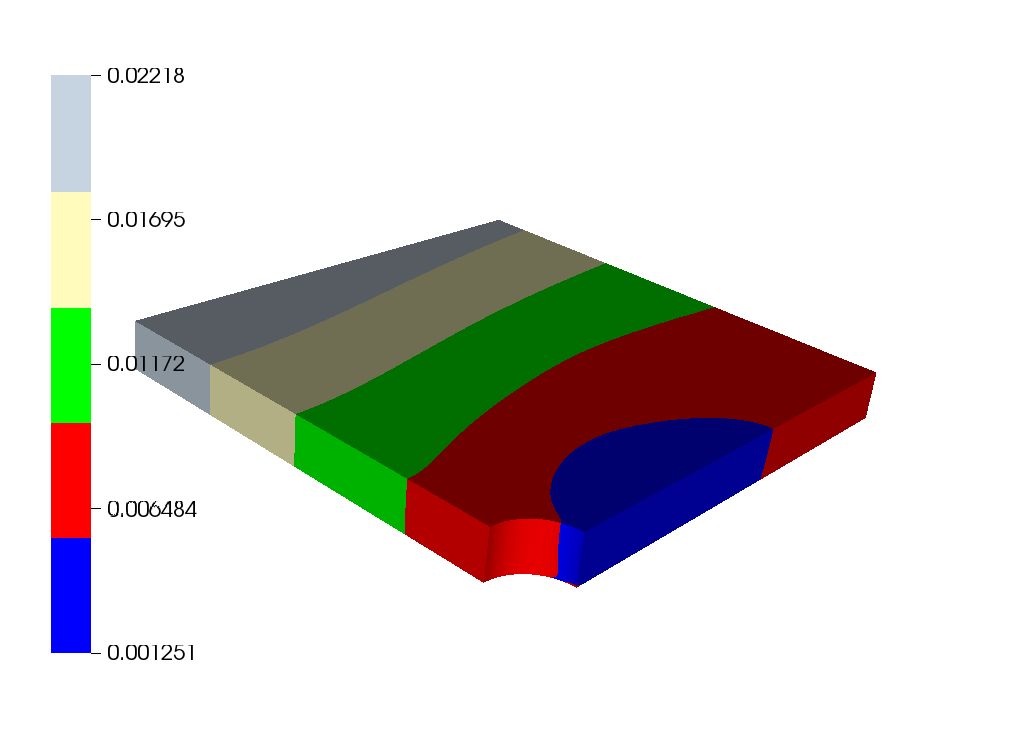}
   \caption{\small{Total displacement $\|u\|$ at $t_8$}}
   \label{fig16}
\end{minipage}
\end{figure}

In Table \ref{tab:time_steps}, the computational history of Algorithm \ref{TFETI_elasto_plastic_algorithm} for the time discretization $(b)$ is documented. We consider the finest mesh level, the stopping tolerances (\ref{tolerances}) and the lumped preconditioner only.  The corresponding development of the plastic zone is depicted for the times $t_2$, $t_4$, $t_6$, and $t_8$ in Figures \ref{fig11} - \ref{fig14}, respectively. The growing zone of plastic elements results from the monotonically increased loading.
Distributions of the von Mises stress $\|\mbox{dev}(\sigma)\|_F$ and the total displacement $\|u\|$ at the final time $t^*$ are depicted in Figures \ref{fig15} and \ref{fig16}, respectively.


\begin{table}[n] 
\begin{center}
\begin{small}
\begin{tabular}{|l|r|r|r|r|r|}
\hline
\multicolumn{6}{|c|}{Mesh properties} \\ \hline
Mesh level  & 1 & 2 & 3 & 4 & 5\\
Number of mesh nodes             &   520 &  1 623 &  9 947 & 166 374 &   309 546 \\
Number of mesh elements          & 1 441 &  6 279 & 48 287 & 931 709 & 1 758 907 \\
\hline
\multicolumn{6}{|c|}{Domain Decompositions} \\ \hline
Number of subdomains             &     1 &      1 &      8 &     135 &       270 \\
Number of CPU cores              &     2 &      2 &      9 &      28 &        28 \\
Number of primal variables       & 1 560 &  4 869 & 33 414 & 623 904 & 1 183 101 \\
Number of dual variables         &   284 &    618 &  5 448 & 136 904 &   272 053 \\
\hline
\multicolumn{6}{|c|}{Performace using the lumped preconditioner} \\ \hline
Number of plastic elements       &   252 &  2 152 & 20 790 & 420 760 &   796 520 \\
Number of Newton iterations.     &     5 &      6 &      6 &       6 &         6 \\
Total number of PCGP iterations  &   237 &    391 &    573 &     718 &       724 \\
Total time in seconds            &     8 &     60 &     77 &     928 &      1796 \\
\hline
\multicolumn{6}{|c|}{Performace using the Dirichlet preconditioner} \\ \hline
Number of plastic elements       &   252 &  2 152 & 20 790 & 420 759 &   796 517 \\
Number of Newton iterations.     &     5 &      6 &      6 &       6 &         6 \\
Total number of PCGP iterations  &   235 &    339 &    440 &     482 &       486 \\
Total time in seconds            &     9 &    174 &    207 &    1532 &      3274 \\
\hline
\end{tabular}
\caption{Numerical results of Algorithm \ref{TFETI_elasto_plastic_algorithm} for different mesh levels and the time discretization $(a)$.}
\label{tab:diferent_meshes}
\end{small}
\end{center}
\end{table}




\begin{table}[n] 
\begin{center}
\begin{small}
\begin{tabular}{|c|c|c|r|r|r|}
\hline
Newton iter. & Number of   &  Number of     & Stopping\\
             & PCGP iters. &  plastic elem. & criterion  \\
\hline
1 & 224 &       0 &          1 \\
2 & 312 & 808 116 &  1.4016e-1 \\
3 & 310 & 797 081 &  4.6051e-2 \\
4 & 319 & 788 802 &  6.3207e-3 \\
5 & 318 & 795 245 &  4.5604e-4 \\
6 & 319 & 796 519 &  1.0735e-5 \\
7 & 321 & 796 596 &  1.1790e-8 \\
8 & 314 & 796 596 & 9.0488e-14 \\
\hline
\end{tabular}
\caption{Computational history of Algorithm \ref{TFETI_elasto_plastic_algorithm} for the time discretization $(a)$ and the finest mesh level}
\label{tab:martin5_upr}
\end{small}
\end{center}
\end{table}


\begin{table}[ht]
\begin{center}
\begin{small}
\begin{tabular}{|r|c|c|r|r|}
\hline
Time  & Number of     & Total number of & Number of      & Solution \\
step  & Newton iters. & PCGP iterations & plastic elems. & time\\
\hline
1 & 2 & 226 &       0 &   536 \\
2 & 2 & 228 &   8 700 &   534 \\
3 & 4 & 487 & 181 557 & 1 152 \\
4 & 4 & 522 & 316 718 & 1 237 \\
5 & 4 & 553 & 502 605 & 1 270 \\
6 & 5 & 863 & 671 067 & 1 678 \\
7 & 5 & 899 & 762 501 & 1 638\\
8 & 4 & 662 & 799 989 & 1 445\\
\hline
\end{tabular}
\caption{Computational history of Algorithm \ref{TFETI_elasto_plastic_algorithm} for the time discretization $(b)$ and the finest mesh level}
\label{tab:time_steps}
\end{small}
\end{center}
\end{table}



Comparing time discretizations $(a)$ and $(b)$ we see that the resulting number of plastic elements differs at the final time $t^*$. The difference is less than $0.5\%$ and is caused by the roundoff errors, the use of iterative solvers, and the numerical evaluation of the yield function which decides whether an element plasticizes or not.

\section{Conclusion}
We have proposed the algorithm for the efficient parallel implementation of elastoplastic problems with the von Mises plastic criterion and the linear isotropic hardening law which is based on the TFETI domain decomposition method. For the time discretization we used implicit Euler method and for the space discretization of the respective one time step elastoplastic problem the finite element method. The latter results in a system of nonlinear equations with a strongly semismooth and strongly monotone operator. Thus the semismooth Newton method was applied and respective linearized problems were solved in parallel using TFETI. The performance of our algorithm was demonstrated on the 3D elastoplastic thin plate with the hole in the center and prescribed loading history. Numerical results for different time discretizations and mesh levels were presented and discussed. Local quadratic convergence of the semismooth Newton method was observed after identification of the plastic zone which is in accordance with the theoretical results.

Our future plan is to apply the  proposed TFETI based algorithm on elastoplastic multi-body contact problems of mechanics. This has already been done successfully for the pure elastic case \cite{DosKozVonBrzMarIJNME-2009,DosKozMarBrzVonHor-CMAME-2010}.

\section*{Acknowledgements}
This work was supported by the European Regional Development Fund in the IT4Innovations Centre of Excellence project (CZ.1.05/1.1.00/02.0070), the project SPOMECH reg. no. CZ.1.07/2.3.00/20.0070, and POSTDOCI II reg. no.  CZ.1.07/2.3.00/30.0055 within Operational Programme Education for Competitiveness. The first, the third and the fourth author acknowledge the support of the project 13-18652S (GA CR) and the second author also thanks to the project of major infrastructures for research, development and innovation of Ministry of Education, Youth and Sports with reg. num. LM2011033 and 7E12028. Authors would particularly like to thank to T. Brzobohaty and A. Markopoulos (Ostrava) for advice  concerning the numerical implementation. A particular thank should go to the unknown reviewers, whose suggestions helped to improve the text significantly.

\bibliographystyle{plain}

\begin{thebibliography}{99}

\bibitem[ACFK02]{AlbCarFun-2002} Alberty, J., Carstensen, C., Funken, S.A., Klose, R.,
\textit{Matlab implementation of the finite element method in elasticity}
Computing 69 (3), 239--263,  2002.

\bibitem[ACZ99]{ACZ99}
Alberty, J., Carstensen, C., Zarrabi, D.,
\textit{Adaptive numerical analysis in primal elastoplasticity with hardening},
Comput. Meth. Appl. Mech. Eng. 171 (3--4), 175--204, 1999.

\bibitem[AC02]{AC02}
Alberty, J., Carstensen, C.,
\textit {Discontinuous Galerkin time discretization in elastoplasticity: motivation, numerical algorithms, and applications},
Comput. Meth. Appl. Mech. Eng. 191 (43), 4949--4968, 2002.

\bibitem[BG94]{BG94}
Badea, L., Gilormini, P.,
\textit{Application of a domain decomposition method to elastoplastic problems},
Int. J. Solids Struct. 31 (5), 643--656, 1994.

\bibitem[BGL05]{BenGolLie05} Benzi, M., Golub, G. H., Liesen, J.,
\textit{Numerical solution of saddle point problems},
Acta Numerica, 1--137, Cambridge University Press, 2005.

\bibitem[Bl99]{BLA99} Blaheta, R.,
\textit{Numerical methods in elasto-plasticity},
Documenta Geonica 1998, PERES Publishers, Prague, 1999.

\bibitem[BCV04]{BroCarVal04} Carstensen, C., Brokate, M., Valdman J.,
{\it A quasi-static boundary value problem in multi-surface elastoplasticity. I: Analysis},
Math. Meth. Appl. Sci. 14 (27), 1697--1710, 2004.

\bibitem[BCV05]{BroCarVal05} Carstensen, C., Brokate, M., Valdman, J.,
{\it A quasi-static boundary value problem in multi-surface elastoplasticity. II: Numerical solution},
Math. Meth. Appl. Sci. 8 (28), 881--901, 2005.

\bibitem[BS96]{BS96}
Brokate, M., Sprekels, J.,
\textit{Hysteresis and Phase Transitions},
Springer, 1996.

\bibitem[BDKKM10]{BrzDosKovKozMar-2010} Brzobohat\'y, T., Dost\'al, Z., Kov\'a\v{r}, P., Kozubek, T., Markopoulos, A.,
\textit{Cholesky decomposition with fixing nodes to stable evaluation of a generalized inverse of the stiffness matrix of a floating structure},
Int. J. Numer. Methods Eng. 88 (5), 493--509,  2011.

\bibitem[CK02]{CarKlo-2002} Carstensen, C., Klose, R.,
\textit{Elastoviscoplastic finite element analysis in 100 lines of Matlab}
J. Numer. Math. 10 (3) , 157--192, 2002.

\bibitem[CHS13]{CHS13}
\v{C}erm\'ak, M., Haslinger, J., Sysala, S.,
 \textit{Numerical solutions of perfect plastic problems with contact: PART II - numerical realization}, in  Proceedings of the XII International Conference on Computational Plasticity - Fundamentals and Applications (eds. E. Onate, D.R.J. Owen, D. Peric, B. Suarez), Barcelona, Spain, 999-1009, 2013.

\bibitem[CKM11]{CerKozMar-2011} \v{C}erm\'ak, M., Kozubek, T., Markopoulos, A.,
\textit{An efficient FETI based solver for elasto-plastic problems of mechanics},
Computational Plasticity XI - Fundamentals and Applications, COMPLAS XI, 1330--1341, 2011.


\bibitem[Doh03]{Dohrmann-2003}  Dohrmann, C. R.,
\textit{A preconditioner for substructuring based on constrained energy minimization}, SIAM J. Sci. Comput., 25, 246--258, 2003.


\bibitem[D09]{DosOQPA-2007} Dost\'al, Z.,
\textit{Optimal Quadratic Programming Algorithms, with Applications to Variational Inequalities},
SOIA 23, Springer US, 2009.

\bibitem[DHK06]{TFETI2006} Dost\'al, Z., Hor\'ak, D., Ku\v{c}era, R.,
\textit{Total FETI - an easier implementable variant of the FETI method for numerical solution of elliptic PDE},
Commun. Numer. Methods Eng. 22 (12),  1155--1162, 2006.

\bibitem[DKMBVH12]{DosKozMarBrzVonHor-CMAME-2010} Dost\'al, Z., Kozubek, T., Markopoulos, A., Brzobohat\'y, T., Vondr\'ak, V., Horyl, P.,
\textit{Theoretically supported scalable TFETI algorithm for the solution of multibody 3D contact problems with friction},
Comput. Meth. Appl. Mech. Eng. 205, 110--120, 2012.

\bibitem[DKVBM10]{DosKozVonBrzMarIJNME-2009} Dost\'al, Z., Kozubek, T., Vondr\'ak, V., Brzobohat\'y, T., Markopoulos, A.,
\textit{Scalable TFETI algorithm for the solution of multibody contact problems of elasticity},
Int. J. Numer. Methods Eng. 82, 1384--1405, 2010.

\bibitem[FR91]{fr91} Farhat, C., Roux, F-X.,
\textit{A method of finite element tearing and interconnecting and its parallel solution algorithm},
Int. J. Numer. Methods Eng. 32, 1205--1227, 1991.


\bibitem[FR92]{fr92} Farhat, C., Roux, F-X.,
\textit{An unconventional domain decomposition method for an efficient parallel solution of large-scale finite element systems},
SIAM J. Sci. Stat. Comput. 13, 379--396, 1992.

\bibitem[FMR94]{fmr94} Farhat, C., Mandel, J.,  Roux, F-X.,
\textit{Optimal convergence properties of the FETI domain decomposition method},
Comput. Meth. Appl. Mech. Eng. 115, 365--385, 1994.


\bibitem[FLLPR01]{FLLPR01}
Farhat, C., Lesoinne, M., LeTallec, P., Pierson, K., Rixen, D.,
\textit{FETI-DP: A dual-primal unified FETI method, part i: A faster alternative to the two-level FETI method},
Int. J. Numer. Methods Eng.  50, 1523--1544, 2001.

\bibitem[FLR00]{FarLesPie-2000} Farhat, C., Lesoinne, M., Pierson, K.,
\textit{A scalable dual-primal domain decomposition method},
Numer. Linear Algebr. Appl. 7, 687--714, 2000.

\bibitem[FP03]{FraPap-2003}
Fragakis, Y., Papadrakakis, M.,
\textit{The mosaic of high performance domain Decomposition Methods for Structural Mechanics: Formulation, interrelation and numerical efficiency of primal and dual methods},
Comput. Methods Appl. Mech. Engrg. 192, 3799--3830, 2003.

\bibitem[FK80]{FukKuf80} Fu\v{c}{\'i}k, S., Kufner, A.,
\textit{Nonlinear Differential Equation}, Elsevier, 1980.

\bibitem[GKLSV12]{GKLSV12}Gruber, P., Kienesberger, J., Langer, U., Sch\"oeberl, J., Valdman, J.,
\textit{Fast solvers and a posteriori error estimates in elastoplasticity},
Langer, Ulrich (ed.) et al., Numerical and symbolic scientific computing. Progress and prospects,
Springer. Texts \& Monographs in Symbolic Computation, 45--63, 2012.

\bibitem[GV09]{GV09} Gruber, P.,  Valdman, J.,
\textit{Solution of One-Time Step Problems in Elastoplasticity by a Slant Newton Method}, SIAM J. Sci. Comp. 31, 1558--1580, 2009.

\bibitem[HR99]{HaRe99} Han, W., Reddy B. D.,
\textit{Plasticity: mathematical theory and numerical analysis}, Springer, 1999.

\bibitem[HV07]{HofVal07} Hofinger, A., Valdman, J.,
\textit{Numerical solution of the two-yield elastoplastic minimization problem},
Computing 81, No. 1, 35--52, 2007.



\bibitem[Jo76]{J76} Johnson, C.,
\textit{Existence theorems for plasticity problems},
J. Math. pures et appl. 55, 431--444, 1976.

\bibitem[Jo78]{J78} Johnson, C.,
\textit{On plasticity with hardening},
J. Math. Anal. Appl. 62, 325--336, 1978.

\bibitem[JPF97]{JPF97} Justino JR, M.R., Park, K.C., Felippa, C.A.,
\textit{An algebraically partitioned FETI method for parallel structural analysis: Implementation and numerical performance evaluation},
Int. J. Numer. Meth. Engng., 40, 2739--2758, 1997.

\bibitem[KLV04]{KLV04}
Kienesberger J., Langer U., Valdman J.,
\textit{On a robust multigrid-preconditioned solver for incremental plasticity problems},
Proceedings of IMET 2004 - Iterative Methods, Preconditioning \& Numerical PDEs, Prague.


\bibitem[KR06]{KlaRhe-2006} Klawonn, A., Rheinbach, O.,
\textit{A parallel implementation of Dual-Primal FETI methods for three dimensional linear elasticity using a transformation of basis},
SIAM J. Sci. Comput. 28, 1886--1906, 2006.

\bibitem[KR10]{KlaRhe-2010} Klawonn, A., Rheinbach, O.,
\textit{Highly scalable parallel domain decomposition methods with an application to biomechanics}, Z. Angew. Math. Mech. 90 (1), 5--32, 2010.

\bibitem[KW01a]{KW01a} Klawonn, A., Widlund, O. B.,
\textit{FETI and Neumann-Neumann Iterative Substructuring Methods: Connections and New Results},
Communications on Pure and Applied Mathematics 54, 57--90, 2001.

\bibitem[KW01b]{KlaWid-2001} Klawonn, A., Widlund, O. B.,
\textit{Dual and dual-primal FETI methods for elliptic problems with discontinuous coefficients in three dimensions}, in: Domain Decomposition Methods, Proceedings of the 12th International Conference on Domain Decomposition
Methods, Chiba, Japan, October 1999, (DDM.org, Augsburg, Germany, 2001).

\bibitem[KL84]{KL84}
Korneev, V. G., Langer U.,
\textit{Approximate solution of plastic flow theory problems},
Teubner-Verlag, Leipzig, volume 69, 1984.

\bibitem[KMBKVD09]{MatSol-2009} Kozubek, T., Markopoulos, A., Brzobohat\'y, T., Ku\v{c}era, R., Vondr\'ak, V., Dost\'al, Z.
\textit{MatSol - MATLAB efficient solvers for problems in engineering},
http://matsol.vsb.cz/.

\bibitem[KVMHDHKC12]{KVMHDHKC} Kozubek, T., Vondr\'ak, V., Men\v{s}\'ik M., Hor\'ak D., Dost\'al, Z., Hapla V., Kabel\'ikov\'a P., \v{C}erm\'ak M.,
\textit{Total FETI domain decomposition method and its massively parallel implementation}, Adv. Eng. Softw. 60-61, 14--22, 2013.

\bibitem[Kr96]{K96}
Krej\v c\'i, P.,
\textit{ Hysteresis, Convexity and Dissipation in Hyperbolic Equations},
GAKUTO Inter. Ser., Math. Scie. Appls., 1996.

\bibitem[KKM13]{KucKozMar-2013}
Ku\v{c}era, R., Kozubek, T., Markopoulos, A.,
\textit{On large-scale generalized inverses in solving two-by-two block linear systems},
Linear Alg. Appl., 438 (7), pp. 3011-3029, 2013.

\bibitem[Man93]{Mandel-1993}
Mandel, J.,
\textit{Balancing domain decomposition},
Commun. Numer. Methods Eng., 9, 233--241, 1993.

\bibitem[MD03]{ManDoh-2003}
Mandel, J., Dohrmann, C. R.,
\textit{Convergence of a balancing domain decomposition by constraints and energy minimization},
Numer. Linear Algebr. Appl. 10, 639--659, 2003.

\bibitem[MT96]{ManTez-1996}
Mandel, J., Tezaur, R.,
\textit{Convergence of a substructuring method with Lagrange multipliers},
 Numer. Math. 73, 473--487, 1996.



\bibitem[Ma79]{Ma79}
Matthies, H.,
\textit{Existence theorems in thermoplasticity},
J. Mecanique 18, No. 4, 695--712,  1979.

\bibitem[Ma79b]{Ma79b}
Matthies, H.,
\textit{Finite element approximations in thermo-plasticity},
Numer. Funct. Anal. Optim. 1, No. 2, 145--160,  1979.

\bibitem[Mi77]{Mi77}
Mifflin, R.,
\textit{Semismoothness and semiconvex function in constraint optimization},
SIAM J. Cont. Optim. 15, 957--972, 1977.

\bibitem[MUMPS]{MUMPS}
“MUMPS Web page”, http://graal.ens-lyon.fr/MUMPS/

\bibitem[NH81]{NeHl81}
Ne\v{c}as, J., Hlav\'{a}\v{c}ek, I.,
\textit{Mathematical Theory of Elastic and Elasto-Plastic Bodies. An Introduction}, Elsevier,  1981.

\bibitem[PJF97]{PJF97} Park, K.C., Justino Jr., M.R., Felippa, C.A.,
\textit{An algebraically partitioned FETI method for parallel structural analysis: Algorithm description},
Int. J. Numer. Meth. Engng., 40, 2717--2737, 1997.

\bibitem[QS93]{QS93}  Qi, L., Sun, J.,
\textit{A nonsmooth version of Newton' s method}, Math. Progr., 58, 353--367, 1993.

\bibitem[RF99]{RF99} Rixen, D.J., Farhat, C.,
\textit{A simple and efficient extension of a class of substructure based preconditioners to heterogeneous structural mechanics problems},
Int. J. Numer. Methods Eng.  44, 489--516, 1999.

\bibitem[Ro92]{r92} Roux, F-X.,
\textit{Spectral analysis of interface operator},
Proceedings of the 5th Int. Symp. on Domain Decomposition Methods for Partial
Differential Equations, ed. D. E. Keyes et al., SIAM, 73--90, Philadelphia 1992.

\bibitem[RF98]{r98} Roux, F-X., Farhat, C.,
\textit{Parallel implementation of direct solution strategies for the coarse grid solvers in 2-level FETI method},
Contemporary Math. 218, 158--173, 1998.

\bibitem[SW11]{SW11}
Sauter, M., Wieners, C.,
\textit{On the superlinear convergence in computational elasto-plasticity},
Comput. Meth. Appl. Mech. Eng. 200, 3646--3658, 2011.

\bibitem[SH98]{SimoHughes}
Simo, J. C., Hughes, T. J. R.,
\textit{Computational Inelasticity},
Interdiscip. Appl. Math. 7, Springer-Verlag, 1998.

\bibitem[NPO08]{NPO08}
de Souza Neto, E. A., Peri\'c, D., Owen, D. R. J.,
\textit{Computational methods for plasticity: theory and application.} Wiley, 2008.

\bibitem[St03]{Stein}
Stein, E. et al.,
\textit{Error-Controlled Adaptive Finite Elements in Solid Mechanics},
Wiley, 2003.

\bibitem[Sy12]{Sy09}
Sysala, S.,
\textit{Application of a modified semismooth Newton method to some elasto-plastic problems}, Math. Comput. Simul. 82, 2004--2021, 2012.

\bibitem[Sy13]{Sy12}
Sysala, S.,
\textit{Properties and simplifications of constitutive time-discretized elastoplastic operators}, Z. Angew. Math. Mech., 1--23, 2013, DOI: 10.1002/zamm.201200056.






\bibitem[TW05]{TosWid-2005}Toselli, A., Widlund, O. B.,
\textit{Domain Decomposition Methods - Algorithms and Theory}, Springer, 2005. 

\bibitem[Wi10]{Wi10}
Wieners, C.,
\textit{A geometric data structure for parallel finite elements and the application to multigrid methods with block smoothing},
Computing and Visualization in Science 13 (4), 161--175, 2010.

\end{thebibliography}

\newpage
\listoffigures 
\listoftables

\end{document}